\documentclass[12pt]{article}
\usepackage{graphicx} % Required for inserting images
\usepackage{amsmath,amssymb,amsthm}
\usepackage{fullpage}

\usepackage{todonotes}

\usepackage{hyperref}

\usepackage[mathlines]{lineno}
\newcommand{\cL}{\mathcal{L}}

\newtheorem{theorem}{Theorem}[section]
\newtheorem{lemma}{Lemma}[section]
\newtheorem{corollary}[theorem]{Corollary}
\DeclareMathOperator{\sech}{sech}

\newcommand{\R}{\mathbb{R}}

\title{The role of coupling and timescales for interacting tipping elements}

\author{Paul D. L. Ritchie\thanks{Department of Mathematics and Statistics, University of Exeter, Exeter EX4 4QF, UK} \and Robbin Bastiaansen\thanks{Institute for Marine and Atmospheric research Utrecht, Department of Physics, Utrecht University, Netherlands} \thanks{Mathematical Institute, Utrecht University, Netherlands} \thanks{Centre for Complex Systems Studies, Utrecht University, Utrecht, The Netherlands} \and Anna S. von der Heydt${}^{\dag}{}^\mathsection$ \and Peter Ashwin${}^{*}$}
\date{July 2026}

\begin{document}

\maketitle

% Other possible titles
% \begin{itemize}
% \item Conditions for tipping propagation in a coupled system
% \item Understanding nonlinear tipping behaviour in a coupled system
% \item Accelerating tipping cascades and the influence of coupling
% \item Propagation of tipping in an accelerating cascade
% \item The role of coupling in interacting tipping elements
% \end{itemize}

\begin{abstract}
Sudden and abrupt changes can occur in a nonlinear system within many fields of science when such a system crosses a tipping point; rapid changes of the system can then occur in response to slow changes in an external forcing. These can occur when time-varying inputs cross a bifurcation. If an ``upstream'' system loses stability in this way, it may cause a ``downstream'' system influenced by it to tip. This especially happens if the downstream system evolves on a much faster timescale, in what we call an accelerating cascade of tipping elements. In this paper, we identify the conditions on coupling and timescales of these systems that result in such tipping cascades and suggest a taxonomy to classify the different possible timings of tipping sequences. We also present a prototypical example of a unidirectionally coupled pair of simple tipping elements with hysteresis. This allows us to map out the various types of response as a function of system parameters and to link it to bifurcations of the underlying system that may have multiple timescales.
\end{abstract}

\newpage 

\tableofcontents

\section{Introduction}

Tipping can occur when a nonlinear system is subjected to slow changes in an external forcing beyond a critical threshold \cite{Lenton2008,armstrong2022exceeding,lenton2023global}. At the critical threshold, the stable base state the system is tracking may cease to exist or become unstable and therefore cause an abrupt transition to an alternative stable state \cite{Ashwin2012tipTypes}. Tipping phenomena are prominent in many fields of science, including climate \cite{rahmstorf2024atlantic,flores2024critical,petrini2025topographically}, ecological \cite{hessen2023lake,basak2024partial,pearce2025considerations}, and social \cite{schulze2024move,alkemade2024social,lenton2023global} systems. Therefore, understanding possible mechanisms of transmission is vital to prevent dangerous tipping events or to enable positive tipping events.  

The focus on the question of whether a tipping threshold is crossed (or not) can lead the reader to suspect that the tipping is instantaneous and irreversible (once a system has tipped, a simple reversal of the forcing does not recover the system back to its original state) \cite{rockstrom2023safe}. In reality, the process of tipping from one state to another takes some time\cite{ritchie2024esd}, and becomes irreversible on longer time scales \cite{armstrong2022exceeding,wang2023mechanisms}.

In complex or heterogeneous systems, several systems with the potential to tip (or ``tipping elements'') may be linked together and subjected to time-varying inputs, and these tipping elements may have vastly different timescales \cite{wang2023mechanisms,wunderling2024climate}. For example, in the climate system the timescales of tipping elements vary on the order of years (e.g. coral reefs) to the order of millennia (e.g. large ice sheets) compared to emissions that change on a timescale of decades, and understanding the perspective of this variety of timescales is vital to understand how tipping in one element may influence tipping of another. For example, \cite{sinet2024amoc} considers a simple model for the influence of Greenland ice sheet tipping on the Atlantic Meridional Overturning Circulation (AMOC), including likely timescales, and finds that the AMOC tipping is affected by the shape of the transient while the ice sheet tips. In particular, the AMOC can shut down due to increased freshwater flux in the North Atlantic from increased meltwater associated with a tipping point that melts (on a slower timescale) the Greenland ice sheet. In contrast, it may be desirable to trigger positive tipping in social systems to reach the rapid levels of decarbonisation required to tackle the climate crisis \cite{eker2024cross,nijsse2025leverage,lenton2022operationalising,meldrum2023breakthrough}.

Very different behaviour can be observed if the system timescale is relatively slow compared to the forcing timescale - indeed, rapidly changing inputs cause rate-induced effects that can force a system to move far from the attractors of the autonomous or ``frozen system'' \cite{Ritchie2023RateInducedTipping}. One effect is the appearance of ``R-tipping'' even if there are no bifurcations of the frozen system \cite{Ashwin2012tipTypes,ashwin2017parameter}. Another effect is that it may be possible to overshoot a threshold and not initiate tipping \cite{ritchie2021overshooting,bochow2023overshooting}. For this, the maximum peak overshoot of the threshold is approximated to be inversely proportional to the square of the overshoot duration, in cases where the forcing timescale is more rapid than the system response \cite{ritchie2019inverse}.

Similarly, different behaviour can be observed based on the relative timescales of the connected systems. Tipping in one system causing another system's threshold to be crossed and hence to tip has been referred to as a domino cascade \cite{Klose2021Cascade}. However, these ``dominoes'' could fall at different speeds: if the first system has a faster timescale than the second, then the first domino would fall before the second, but if these timescales are reversed, then the second domino will topple much quicker than the first. %Understanding and developing a taxonomy of the different tipping sequences that can arise is one key focus of this paper.
%{\color{red}RB: I'd add here a paragraph to talk specifically about the relevance of time scales and timing etc. for tipping cascades. So what could happen and that during the tiping something can happen to the other system etc. That seems what ref 2 wants, and I think one paragraph could be enough. } {\color{blue}PR: I've had a go at this above, referencing Klose et al., 2021, but to a large extent I reckon we are doing this in 2 paragraphs (penultimate paragraph of intro) down from here...?}

The current paper introduces a taxonomy for {\em tipping sequences}, i.e. possible sequences of tipping events of an upstream and a downstream system.
%, depending on the functional form of the coupling as well as the system and forcing timescales. 
This taxonomy recognizes that tipping onset time and duration depend on system timescale, coupling and forcing. We use this to identify the conditions for different sequences of upstream and downstream tipping that can arise from a single monotonic external forcing profile applied to the upstream system. We illustrate this taxonomy using a deterministic version of the model in \cite{ashwin2025early} which is studied here in a lot more depth.

In \cite{ashwin2025early} we explore what it means for an early warning signal to give a skilful prediction of tipping using an extrapolation methodology for critical slowing-down based indicators, applied to an idealised cascade with only two coupled tipping elements modelled as a non-autonomous stochastic differential equation. The model consists of an (upstream) subsystem that is forced by time-varying input, which its output forces another (downstream) subsystem.
In particular, we show in \cite{ashwin2025early} how it is difficult to get skilful warning of ``downstream tipping within upstream tipping'' because of breakdown of extrapolation of downstream indicators as the upstream system starts to tip.

The rest of this paper is set up as follows. Section~\ref{sec:timescales_cascades} explains the scenario in general terms, while Section~\ref{sec:propagation} gives precise criteria for the upstream and downstream systems to undergo tipping. This helps explain the influence of coupling on the tipping of upstream and downstream systems. In particular, we classify the temporal scenarios for upstream and downstream tipping, noting that the duration of tipping is finite in both cases, and may be different depending on the relative timescales of the systems. Section~\ref{sec:example} looks at a specific example where we classify the tipping regimes of coupled one-dimensional systems based on the coupling between systems and timescale separation. We show how the dynamics change based on the type of coupling between the subsystems (including the timing of tipping relative to each other), examining both linear and localised state coupling that gives rise to overshoots of the downstream system's threshold. To illustrate this in detail, Section~\ref{sec:linear} explores this for linear coupling, and Section~\ref{sec:localised} for a form of localised coupling. We find that the strength and timescale separation (between the upstream and downstream systems) can give rise to many different tipping (cascade) regimes, even in the linear case. Finally, we explore some generalizations and implications of this to other systems and possible applications in Section~\ref{sec:discuss}.

\subsection{Timescales and cascades of tipping elements}\label{sec:timescales_cascades}

\begin{table}[t]
	\centering
    \caption{\textbf{Summary of notation.} We summarise the conventions used in the notation for systems such as the coupled systems in Figure~\ref{fig:block_diag}.
}
	\begin{tabular}{l|p{13cm}}
	\hline
	$x$ & state of upstream system\\
	$y$ & state of downstream system\\
	$t$ & rescaled time coordinate (taken in reference to the upstream system)\\
    $t'$ & rescaled time coordinate (taken in reference to the downstream system) \\
	$\Lambda$ & forcing function for upstream system\\
	$f$ & function determining upstream system dynamics\\
	$M$ & coupling function from upstream to downstream system\\
	$g$ & function determining downstream system dynamics\\
	$r$ & ratio of timescales: forcing to upstream system\\
	$\epsilon$ & ratio of timescales: upstream to downstream systems\\
    $\lambda_-$ & past limit of the forcing, i.e. $\lambda_- := \lim_{s \rightarrow -\infty} \Lambda(s)$\\
    $\lambda_+$ & future limit of the forcing, i.e. $\lambda_+ := \lim_{s \rightarrow -\infty} \Lambda(s)$ \\
    $\lambda_u$ & location of a fold bifurcation in the examples in the paper \\
    $X_-$ & past limit of the upstream system \eqref{eq:f-nonaut} \\
    $x^{[r]}$ & $r$-dependent pullback attractor of the upstream system \eqref{eq:f-nonaut} limiting to $X_-$ for $t \rightarrow -\infty$\\
    $X_+^{[r]}$ & $r$-dependent future limit of the pullback attractor of the upstream system, i.e. $X_+^{[r]}:=\lim_{t \rightarrow \infty} x^{[r]}(t)$\\
    $(x^{[r]},y^{[r,\epsilon]})$ & $r$ and $\epsilon$-dependent pullback attractor of the coupled system \eqref{eq:slowfastode2} limiting to $(X_-,Y_-)$ for $t \rightarrow -\infty$\\
    $(X_-,Y_-)$ & assumed past limit of the pullback attractor $(x^{[r]},y^{[r]})$ for the coupled system \eqref{eq:slowfastode2}\\
    $\tilde{x}$ & frozen tipping trajectory of the upstream system at the fold bifurcation at $\lambda_f$.\\
	\hline
	\end{tabular}

\label{tab:notation}
\end{table}

In this subsection, we first briefly motivate the setup we consider in this paper where there is forcing of an upstream system that in turn forces a downstream system. Throughout this paper we will consistently denote the upstream system by $x$ and the downstream system by $y$ -- see also Table~\ref{tab:notation} for further notational conventions used throughout the paper. In Section~\ref{sec:propagation}, we will make this more precise and explicit. In this paper, we focus on tipping associated with fold bifurcations and forcing that is asymptotically constant.

\begin{figure}[!ht]
	\centering
	\includegraphics[width=0.6\linewidth]{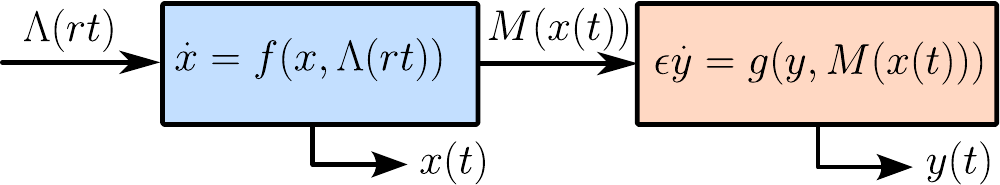}
	
	\caption{\textbf{Block diagram of coupled system.} The slowly varying input $\Lambda(rt)$ forces an upstream system with state variable $x(t)$. The upstream system is coupled to a downstream system with state variable $y(t)$ via the coupling function $M(x(t))$. Parameter $\epsilon>0$ governs the timescale separation between the upstream and downstream systems. In \cite[Fig. 8]{ashwin2025early} we refer to this cascade of tipping elements as accelerating if $\epsilon\ll 1$ and decelerating if $\epsilon \gg 1$.}
	\label{fig:block_diag}
\end{figure}

Specifically, we consider a non-autonomous coupled system whose upstream state $x \in \R^n$ and downstream state $y \in \R^m$ are governed by
\begin{equation}
\left.
\begin{aligned}
    \label{eq:f-nonaut}
   \tau_x\frac{\mathrm{d}x}{\mathrm{d}t'} &= f(x,\Lambda(t'/\tau_\Lambda)) \\
    \tau_y\frac{\mathrm{d}y}{\mathrm{d}t'} &= g(y,M(x(t'/\tau_x)))
\end{aligned}\right\}
\end{equation}
where $t'$ is time, $f:\R^{n+p}\rightarrow \R^n$ and $g:\R^{m+q}\rightarrow \R^q$ are smooth functions. We assume the system is externally forced by $\Lambda:\R\rightarrow \R^p$ at some rate of forcing. A graphical representation of the coupled upstream and downstream system is presented in Figure~\ref{fig:block_diag}.

In this system, there are three time scales at play: (i) the internal time scale $\tau_x$ of the upstream system, (ii) the internal time scale $\tau_y$ of the downstream system and (iii) the time scale of external change $\tau_\Lambda$. As will be discussed in the rest of this paper, the ratio between these time scales influences the cascading tipping dynamics of the full coupled system. To make these ratios more explicit, we will not work with \eqref{eq:f-nonaut} directly, but instead rescale time to the internal time scale of the upstream system by introducing $t := t'/\tau_x$ to obtain the system
\begin{equation}
\left.
\begin{aligned}
    \frac{\mathrm{d}x}{\mathrm{d}t} & = f(x,\Lambda(rt))\\
    \epsilon\frac{\mathrm{d}y}{\mathrm{d}t} & = g(y,M(x(t)))
\end{aligned}\right\}
\label{eq:slowfastode}
\end{equation}
where $r := \tau_x/\tau_\Lambda$ is the ratio between internal time scale of the upstream system and the external forcing, and $\epsilon := \tau_y / \tau_x$ is the ratio between internal time scales of the upstream and downstream systems.

We assume for convenience that $\Lambda(s)$ is asymptotically constant, i.e. that it limits to $\lambda_{\pm}$ as $s\rightarrow \pm\infty$. Throughout this paper, we typically assume the forcing is slow compared to the system dynamics. This corresponds to assuming $r$ is small, as $\tau_x\ll \tau_\Lambda$. 

For such choices of $0<r\ll 1$, we will try to understand the tipping dynamics of the up- and downstream systems of the full coupled system \eqref{eq:f-nonaut} by connecting it to dynamics of autonomous {\em frozen systems} that can be understood by traditional methods.

% upstream system
First, we start by looking at tipping of the upstream system. This tipping of the upstream system in \eqref{eq:slowfastode} can be understood using methods discussed in \cite{Ashwin2012tipTypes,ashwin2017parameter,wieczorek2023rate}; i.e. tipping can be understood in terms of bifurcations of the associated {\em frozen system}
\begin{equation}
\label{eq:f-frozen}
    \frac{\mathrm{d}x}{\mathrm{d}t} = f(x,\lambda)
\end{equation}
on varying parameters $\lambda\in\R^p$.

The idea -- that will be formalised in Section~\ref{sec:propagation} -- is that for varying forcing $\Lambda$ that is slow enough (i.e. $r$ small enough), the non-autonomous system
\begin{equation}
	\label{eq:f-nonaut_tscale}
	\frac{\mathrm{d}x}{\mathrm{d}t} = f(x,\Lambda(rt))
\end{equation}
will behave qualitatively mostly similar to the frozen system \eqref{eq:f-frozen}. Specifically, if $\Lambda$ remains in a region where there is a continuous branch of stable equilibria for (\ref{eq:f-frozen}) then for small enough $r$, solutions will {\em track} the branch of equilibria \cite{ashwin2017parameter}. If on the other hand, $\Lambda$ crosses a fold bifurcation\footnote{Note that tipping associated with other bifurcations is possible, but we do not consider these here \cite{Kuehn2011CSD}.} where a stable equilibrium ceases to exist then solutions will undergo a {\em bifurcation-induced tipping (B-tipping)} for small enough $r$.

This can be made rigorous by application of \cite[Theorem 2.2]{ashwin2017parameter}. For this, we assume that there is a linearly stable equilibrium $X_-$ of the frozen system (\ref{eq:f-frozen}) at $\lambda=\lambda_-$. Then this theorem indicates that there is a solution $x^{[r]}(t)$ of (\ref{eq:f-nonaut_tscale}) that limits to $X_-$ in the past, i.e. $\lim_{t \rightarrow -\infty} x^{[r]}(t) = X_-$. Such solution $x^{[r]}(t)$ is called a pullback attractor, which is a generalization of the notion of attractor in the case of a non-autonomous system, and it describes the (time-dependent) set of points in state space that will be visited for trajectories started at ``typical'' initial conditions in the far past: see for example \cite{KloedenRasmussen2011}. Thus by understanding the (future) dynamics of the pullback attractor $x^{[r]}(t)$, one automatically understands the properties of any solution that starts in the past close to $X_-$.

Hence, we are in particular interested in the future limit of the pullback attractor $x^{[r]}$, which we define as the $r$-dependent state
\begin{equation}
X^{[r]}_+:=\lim_{t\rightarrow\infty} x^{[r]}(t).
\end{equation}
By identifying whether $X^{[r]}$ corresponds to a tipped upstream system or not, we distuingish between tipping or tracking for the upstream system. Specifically, we say there is {\em end-point tracking} of the attractor along a branch of equilibria from $X_-$ to some $X_+=X_+^{[r]}$ for some $r>0$ if there is a branch of linearly stable equilibria $X(\Lambda(s))$ such that 
\begin{equation*}
X_-=X(\lambda_-),~~X_+=X(\lambda_+),
\end{equation*}
and we say there is {\em upstream tipping} if this is not the case.

% downstream system
Second, we inspect the downstream system. For a given solution $x(t)$ of the upstream system, the downstream system is described by the non-autonomous system
\begin{equation}
\epsilon\frac{\mathrm{d}y}{\mathrm{d}t} = g(y,M(x(t))).
\label{eq:g-nonaut_tscale}
\end{equation}
For convenience, we write the downstream system in terms of the timescale $\tilde{t}=t/\epsilon$:
\begin{equation}
\frac{\mathrm{d}y}{\mathrm{d}\tilde{t}} = g(y,M(x(\epsilon \tilde{t})) ),
\label{eq:g-nonaut-tilde}
\end{equation}

We aim to understand the behaviour of (\ref{eq:g-nonaut-tilde}) in relation to the upstream forcing $M(x(\epsilon \tilde{t}))$ and the {\em frozen downstream system}
\begin{equation}
\frac{\mathrm{d}y}{\mathrm{d}\tilde{t}} = g(y,\mu).
\label{eq:g-frozen}
\end{equation}
where $\mu\in\R^q$. The similarity of (\ref{eq:f-nonaut_tscale}) and (\ref{eq:g-nonaut-tilde}) means we can apply similar methods to understand cases where tipping in (\ref{eq:f-nonaut_tscale}) does, or does not lead to tipping in (\ref{eq:g-nonaut-tilde}). In particular, for the downstream system the function $M(x(\epsilon \tilde{t}))$ and the rate $\epsilon$ take on the roles of the function $\Lambda(t)$ and the rate $r$ that were used in the upstream system.

Third, we look at the coupled upstream and downstream system \eqref{eq:slowfastode}. Here, the ideas for the upsteam system hold as presented before. However, for the downstream system, the precise solution $x(t)$ of the upstream system is needed. By the nature of the pullback attractor $x^{[r]}(t)$, and theorems from \cite{Ashwin2012tipTypes,ashwin2017parameter,wieczorek2023rate}, it will turn out that we can use the ideas for the downstream system when we consider it to be forced by the upstream pullback attractor $x^{[r]}$. Since that upstream pullback attractor depends on $r$, the downstream pullback attractor $y^{[r,\epsilon]}$ will as well, in addition to its natural dependence on $\epsilon$.

Now, there are different scenarios depending on the timescale separation $\epsilon$. If $\epsilon \ll 1$ then the downstream dynamics will be much faster than the upstream and we say the cascade {\em accelerates}. If $\epsilon \gg 1$ we say the cascade {\em decelerates}.

Later in Section~\ref{sec:propagation}, for the case of an accelerating cascade of tipping elements ($\epsilon \ll 1$) we find that tipping of the upstream system $x$ may or may not cause tipping of the fast downstream system $y$; Theorem~\ref{thm:criteria} gives criteria for propagation, depending on properties of $M$ and the tipping trajectory of $x$.
For the case of a decelerating cascade ($\epsilon \gg 1$) of tipping elements, it is a moot point whether the downstream system undergoes B-tipping in the usual sense in that, in this case, the changes to forcing $M(x(\epsilon \tau))$ are generally fast in the timescale of the downstream system. If $\epsilon$ is of order one then the upstream and downstream tipping systems have comparably similar timescales and rate-dependent effects can appear. 

\subsection{Timing of downstream relative to upstream tipping}\label{sec:timings}

\begin{figure}[!ht]
    \centering
    \includegraphics[width=0.8\linewidth]{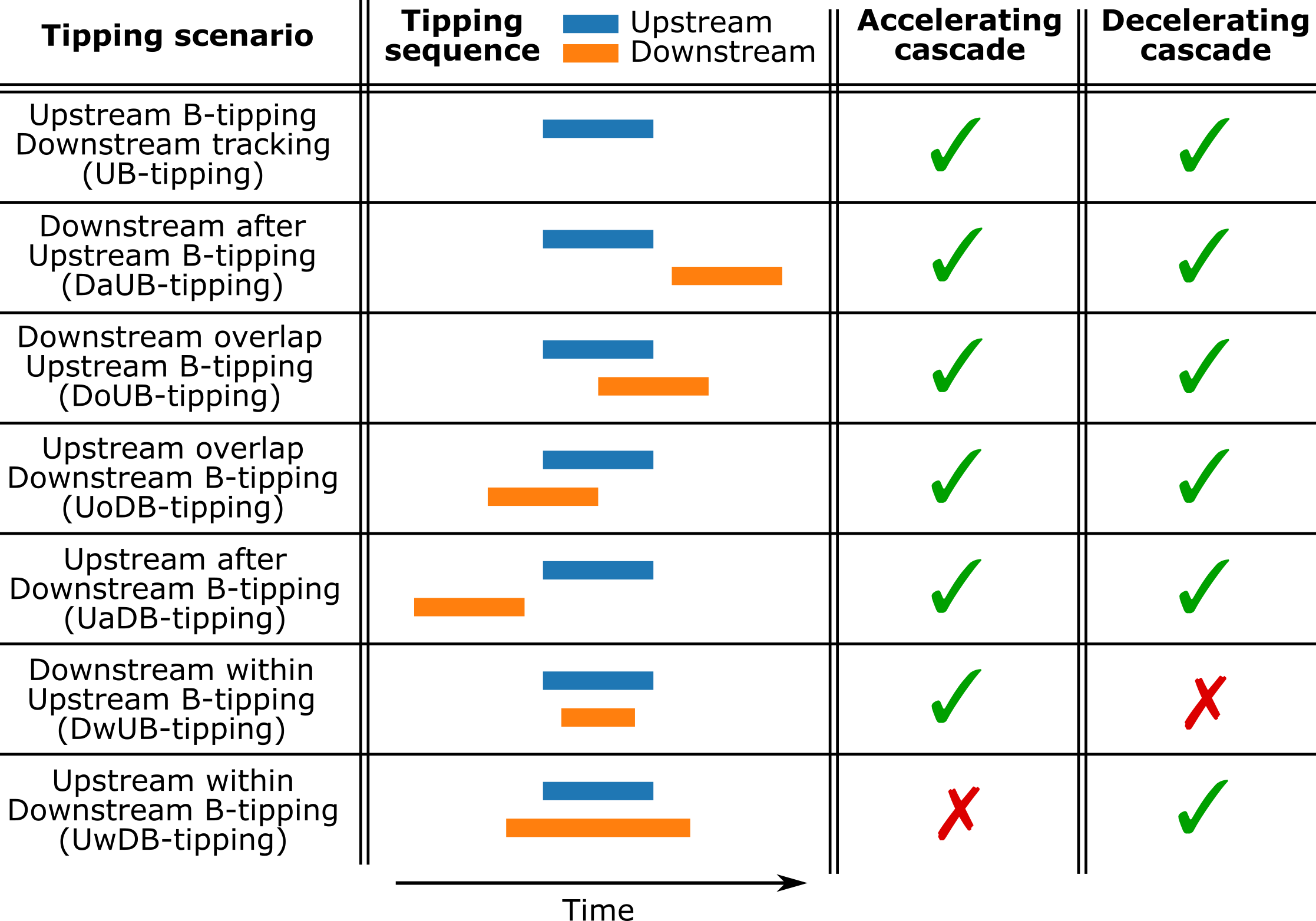}
    
    \caption{\textbf{Schematic of possible tipping sequences for upstream and downstream systems.} The upstream system is assumed to be subjected to the same slowly varying external forcing in all cases, such that tipping of the upstream system occurs at the same time and same duration across all scenarios (blue bar). Tipping of the downstream system (if it occurs, orange bar) changes relative to that of the upstream system tipping. Final two columns indicate if the tipping scenario can occur in accelerating ($\epsilon \ll 1$) and decelerating ($\epsilon \gg 1$) cascades.}
    \label{fig:tipping_seq}
\end{figure}

Many different tipping scenarios can arise based on the timing of the downstream system tipping (if it does) relative to that of the upstream system tipping. These are summarised in Figure~\ref{fig:tipping_seq}, where it is assumed that the external forcing profile is unchanged across scenarios, and such that there is tipping of the upstream system. Note that \cite{ashwin2025early} identified issues that the relative timing can cause for the extrapolation of downstream forcing, and consequently for the predictability of downstream tipping.

We say the tipping starts when the forcing crosses some bifurcation threshold, and the tipping finishes when the state is within some neighbourhood of a new stable branch (we refer to the later Figure~\ref{fig:upstream_cases} for a specific example). For overshooting scenarios where the forcing crosses a threshold but the state does not, we say tipping has not occurred. Clearly, the end of tipping will depend on the choice of neighbourhood. A reasonable size of this neighbourhood could be, for example, chosen with reference to additional stochastic variability that we have not modelled here. %In consequence, the classification is quantitative rather than qualitative. \todo{Other way round?}

Changing the coupling and/or timescale separation between the upstream and downstream dynamics will impact the tipping behaviour of the downstream system relative to the upstream system (though the upstream tipping onset and offset are unaffected, as indicated by the blue bar remaining fixed).

The simplest scenario is when the downstream system does not tip and there is only upstream B-tipping (UB-tipping). If the downstream system does tip, then this can start after the upstream system has completed tipping (DaUB-tipping). Alternatively, the tipping can propagate such that in the process of the upstream system tipping, the downstream system starts tipping but finishes tipping after the upstream system, such that the tipping events overlap (DoUB-tipping). If the downstream system is initially close to its threshold and/or the coupling is strong, then the downstream system could start tipping first. The upstream system could then overlap with the downstream system tipping (UoDB-tipping) or only commence tipping after the downstream system has tipped (UaDB-tipping). The cases described up to this point can occur for either an accelerating or a decelerating cascade. However, unique for an accelerating cascade is that the tipping of the downstream system can be completely contained within the tipping of the upstream system (DwUB-tipping) due to the faster tipping dynamics of the downstream system. On the other hand, unique for a decelerating cascade is that tipping of the upstream system can be contained completely within the downstream tipping (UwDB-tipping). In the following sections, we will provide some examples that demonstrate under what conditions these scenarios emerge. 

\section{Tracking and tipping in cascades of tipping elements}
\label{sec:propagation}

In this Section, we will make the ideas described in Section~\ref{sec:timescales_cascades} rigorous. For this, we consider mechanisms for tipping in a class of systems of the form (\ref{eq:f-nonaut_tscale}, \ref{eq:g-nonaut_tscale}) as shown in Figure~\ref{fig:block_diag}, namely
\begin{equation}
\left.
\begin{aligned}
    \frac{\mathrm{d}x}{\mathrm{d}t} & = f(x,\Lambda(rt))\\
    \epsilon\frac{\mathrm{d}y}{\mathrm{d}t} & = g(y,M(x))
\end{aligned}\right\}
\label{eq:slowfastode2}
\end{equation}
where $x\in\R^n$, $y\in \R^m$, $f:\R^{n+p}\rightarrow\R^n$, $g:\R^{m+q}\rightarrow\R^m$, $\Lambda:\R\rightarrow \R^p$, $M:\R^n\rightarrow \R^q$ are all smooth functions of their arguments.  
For simplicity, we assume that $g(y,\mu)$ does not depend explicitly on $t$ and there is no back-coupling from $y$ to $x$. In the case of $r$ small, this will correspond to quasi-static variation of the parameter in the system $x$. 

We assume $\Lambda(s)$ is a parameter shift, i.e. it is bi-asymptotically constant with
\begin{equation*}
\lambda_{\pm}:=\lim_{s\rightarrow \pm \infty} \Lambda(s)
\end{equation*}
as in \cite{wieczorek2023rate}.

We assume that for $\lambda=\lambda_-$ (i.e. in the limit $t \rightarrow -\infty$) the frozen system
\begin{equation}
\left.
\begin{aligned}
    \frac{\mathrm{d}x}{\mathrm{d}t} & = f(x,\lambda)\\
    \epsilon\frac{\mathrm{d}y}{\mathrm{d}t} & = g(y,M(x))
\end{aligned}\right\}
\label{eq:slowfastodefrozen}
\end{equation}
has a linearly stable equilibrium $(X_-,Y_-)$. Then, by  \cite[Theorem 2.2]{ashwin2017parameter}, for each $r>0$ there is a trajectory $(x^{[r]}(t),y^{[r,\epsilon]}(t))$ of the nonautonomous system (\ref{eq:slowfastode2}) that is a pullback point attractor with limit
\begin{equation*}
\lim_{t\rightarrow -\infty}(x^{[r]}(t),y^{[r,\epsilon]}(t)))=(X_-,Y_-).
\end{equation*}
That is, this trajectory $(x^{[r]}(t),y^{[r,\epsilon]}(t))$ limits in backwards time to the linearly stable equilibrium $(X_-,Y_-)$ of the frozen system with $\lambda = \lambda_-$. In the rest of this section, we aim to understand future limit properties of this trajectory $(x^{[r]}(t),y^{[r,\epsilon]}(t))$ through knowledge of dynamic properties of the frozen system (\ref{eq:slowfastodefrozen}). In particular, this will yield some results that imply downstream tipping of various types described in Section~\ref{sec:timings}.

\subsection{Criteria for tracking and B-tipping}

In this Section, we state some hypotheses and notation that allow us to prove statements about tracking and tipping for upstream and downstream systems (\ref{eq:slowfastode2}). These will depend on the properties of these systems, the forcing and the coupling between them. Although we do not aim to give a complete classification, these methods should be adaptable to a wide range of scenarios. 

We start by considering the upstream system. Here, we focus our attention to attractors that are equilibrium solutions (compared to e.g. periodic orbits), and crossings of fold bifurcations (compared to e.g. Hopf bifurcations). We consider the case that the system at $t \rightarrow -\infty$ was in some linearly stable equilibrium and describe when such system tracks this equilibrium and when it instead tips. Mathematically, we formulate these assumptions as

\begin{itemize}
\item[] {\bf H1 (upstream hypothesis)}
Suppose that all attractors of the frozen upstream system (\ref{eq:f-frozen}) are equilibria and that there is a family of linearly stable equilibrium solutions $X(\lambda)$ for $\lambda$ defined over some connected region $\lambda\in\cL_{X}$. Assume $\lambda_-\in \cL_{X}$ and $X_-:=X(\lambda_-)$, and that $\partial \cL_{X}$ consists of non-degenerate fold bifurcations.
\end{itemize}

Based on the precise forcing profile $\Lambda$ and the forcing time scale $r$, the upstream system can either track the equilibrium (i.e. not tip) or tip to a (potentially) different equilibrium. The general statements we can obtain are given in the following lemma, which is a direct application of \cite[Lemma 2.3]{ashwin2017parameter}, which gives sufficient conditions for behaviour of the pullback attractor $x^{[r]}(t)$ of the upstream system in terms of the forcing $\Lambda$ and rate $r$.

\begin{lemma}[{\bf Sufficient conditions for upstream tracking/B-tipping}]
Suppose that {\bf H1} is satisfied.
\begin{enumerate}
    \item 
    If $\Lambda(s) \subset \cL_{X}$ for all $s$ and $\lambda_+\in \cL_{X}$ then there is an $r_0>0$ such that there is end-point tracking of the upstream system for any $0<r<r_0$.
    \item If there is an $s_0$ such that $\Lambda(s)\not\in \cL_{X}$ for all $s>s_0$ then there is B-tipping of the upstream system for any $r>0$.
    \item 
    If there is an $s_0$ such that $\Lambda(s_0)\not\in \cL_{X}$ but $\lambda_+\in \cL_{X}$ then there is an $r_0>0$ such that there is B-tipping for any $0<r<r_0$.
\end{enumerate}
\label{lem:Btip-up}
\end{lemma}

Case 1 indicates that the upstream system does not tip if the forcing profile $\Lambda$ does not cross a bifurcation \emph{and} the forcing does not change too rapidly. Case 2 implies that if the branch no longer exists beyond some point of time in the forcing, the system cannot remain on it, and so there must be B-tipping for any rate. Finally, case 3 means there is no possibility of a safe overshoot if forcing is slow enough.

We now turn our attention to the downstream system. To be able to meaningfully discuss a cascade of tipping $r$ needs to be chosen such that there is B-tipping of the upstream system -- though we do not restrict ourselves to this case exclusively in the remainder of this section. In order to state sufficient conditions for whether there is also downstream tipping, we need to again make some hypotheses about the downstream system, which are similar to the ones we made for the upstream system:

\begin{itemize}
\item[] {\bf H2 (downstream hypothesis)}
Suppose that all attractors of the frozen downstream system (\ref{eq:g-frozen}) are equilibria and there is a family of linearly stable equilibrium solutions $Y(\mu)$ for $\mu$ defined over some connected region $\mu\in \cL_{Y}$. Assume $\mu_-:=M(X_-)\in \cL_{Y}$ and $Y_-:=Y(\mu_-)$, and that $\partial \cL_{Y}$ consists of non-degenerate fold bifurcations.
\end{itemize}

As with the upstream system, whether the downstream system tips or not depends on the precise forcing profile $M(x^{[r]})$ and the timescale separation parameter $\epsilon$. In particular, if $\epsilon>0$ is small, the question of whether there is propagation of upstream B-tipping to downstream B-tipping depends on whether $M(x^{[r]}(t))$ passes through $\partial \cL_{Y}$. We can again state this as a consequence of \cite[Lemma 2.3]{ashwin2017parameter} by giving sufficient conditions for end-point tracking or B-tipping of the downstream system i.e. whether
\begin{equation*}
\lim_{t\rightarrow \infty} y^{[r,\epsilon]}(t) \rightarrow Y(M(X(\lambda_+))),
\end{equation*}
is satisfied (i.e. end-point tracking), or not (i.e. B-tipping).

\begin{lemma}[{\bf Sufficient conditions for downstream tracking/B-tipping}]
Suppose that {\bf H1} and {\bf H2} are satisfied and fix $r>0$. Then there is a unique pullback attractor $x^{[r]}(t)$ for the upstream system that limits to $X_-$ in the past. Then
\begin{enumerate}
    \item If $M(x^{[r]}(t))\in \cL_{Y}$ for all $t$ and $M(X_+)\in\cL_{Y}$ then there is an $\epsilon_0(r)>0$ such that there is end-point tracking of the downstream system for any $0<\epsilon<\epsilon_0(r)$.
    \item If there is a $t_0$ such that $M(x^{[r]}(t))\not \in \cL_{Y}$ for all $t>t_0$ then there is a B-tipping of the downstream system for any $\epsilon>0$.
    \item If there is a $t_0$ such that $M(x^{[r]}(t_0))\not \in \cL_{Y}$ but $M(X_+)\in \cL_{Y}$ then there is an $\epsilon_0(r) > 0$ such that there is B-tipping for any $0 < \epsilon < \epsilon_0(r)$.
\end{enumerate}
\label{lem:Btip-down}
\end{lemma}

This result states sufficient conditions, respectively, for end-point tracking, for B-tipping (for all $\epsilon$) and for B-tipping in cases where larger $\epsilon$ may lead to rapid overshoot and end-point tracking -- similar to Lemma~\ref{lem:Btip-up} for the upstream system. 

\begin{figure}[!ht]
    \centering
    \includegraphics[width=0.45\linewidth]{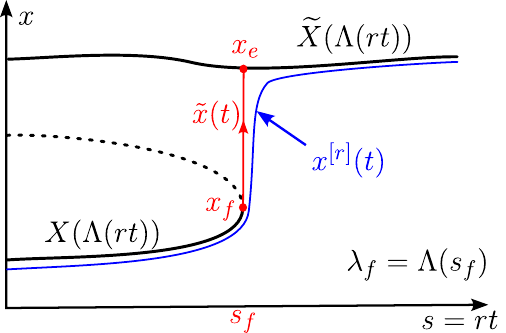}
    
    \caption{\textbf{Approximation of the upstream pullback attractor.} For small $r$, the pullback attractor $x^{[r]}(t)$ of the upstream system (blue line) is approximated by the stable branch $X(\lambda)$ (black) up to the B-tipping where $s_f=rt$. On passing slowly through tipping at $s_f$, the pullback attractor is approximated by the frozen tipping trajectory $\tilde{x}(t)$ (red), which is the connecting trajectory from the saddle node $x_f$ at $\lambda_f=\Lambda(s_f)$ to a stable equilibrium $x_e$ on the new stable branch $\widetilde{X}(\lambda)$ (black). This is used in Theorem~\ref{thm:criteria} to give criteria for the downstream system to undergo B-tipping during upstream B-tipping.}
    \label{fig:schem-proof}
\end{figure}

We now define the frozen tipping trajectory of the upstream system. We will use this specific tipping trajectory of the upstream system to obtain conditions for B-tipping/end-point tracking of the downstream system in the case of slow forcing of an accelerating cascade of tipping elements (i.e. $\epsilon \ll 1$). If the frozen upstream system (\ref{eq:f-frozen}) has a non-degenerate fold bifurcation at $(x_f,\lambda_f)$ where $\lambda_f\in\partial\cL_X$ then there is a unique trajectory $\tilde{x}(t)$ of the upstream frozen system at $\lambda=\lambda_f$ (up to time-translation), such that 
\begin{equation*}
\lim_{t\rightarrow -\infty} \tilde{x}(t)=x_f, ~~ \lim_{t\rightarrow +\infty} \tilde{x}(t)=x_e
\end{equation*}
where $x_e$ is the final attracting equilibrium state of the connecting trajectory for the upstream frozen system. 
We call $\tilde{x}(t)$ a {\em frozen tipping trajectory} of the frozen upstream system (\ref{eq:f-frozen}) at $\lambda=\lambda_f$ (see also Figure~\ref{fig:schem-proof}).

To be able to formulate a meaningful theorem and not further complicate the presentation, we  impose some additional assumptions on the coupled system \eqref{eq:slowfastode2}. In particular, we focus on the crossing of one fold bifurcation in the upstream system, and on the question whether the downstream system tips \emph{during} the tipping of the upstream system. The following hypothesis is constructed such that at most one bifurcation crossing happens in the downstream as well. The following hypothesis has the consequence that, for $r$ and $\epsilon$ small, there is downstream tracking before and after any upstream tipping; in this case the downstream tipping will depend on the effect of coupling via the frozen tipping trajectory only. Specifically, we impose the following:

\begin{itemize}
\item[] {\bf H3 (coupling hypothesis)}
Suppose that {\bf H1} holds, and additionally there is an $s_f$ such that:
\begin{enumerate}
\item For $\Lambda(s_f)=\lambda_f\in\partial \cL_{X}$ there is a frozen tipping trajectory $\tilde{x}$ with $x_f=X(\lambda_f)$.
\item There is a region of linearly stable attracting equilibria $\widetilde{X}(\lambda)$ for $\lambda$ in some connected region $\lambda\in\cL_{\widetilde{X}}$, such that $\Lambda(s)\in \cL_{\widetilde{X}}$ for all $s\geq s_f$.
\end{enumerate} 
Suppose that {\bf H2} holds, and additionally:
\begin{enumerate}
\setcounter{enumi}{2}
\item For all $s<s_f$ we have $M(X(\Lambda(s)))\in\cL_{Y}$.
\item There is a region of linearly stable attracting equilibria $\widetilde{Y}(\mu)$ for $\mu$ in some connected region $\mu\in\cL_{\widetilde{Y}}$, such that $M(\widetilde{X}(\Lambda(s)))\in \cL_{\widetilde{Y}}$ for all $s\geq s_f$.
\item There is a value $t_f$ such that for $\mu_f := M(\tilde{x}(t_f))$ there is a frozen trajectory $\tilde{y}(\tilde{t})$ of the downstream system such that 
\begin{equation*}
\lim_{\tilde{t}\rightarrow-\infty} \tilde{y}(\tilde{t})=Y(\mu_f),~~\lim_{\tilde{t}\rightarrow\infty} \tilde{y}(\tilde{t})=\widetilde{Y}(\mu_f).
\end{equation*}
\end{enumerate} 
\end{itemize}

Note that Hypothesis {\bf H3} 3-5 are assumptions on the coupling from an upstream system undergoing tipping, such that the downstream system can track (possibly different) attractors before and after the upstream tipping. 
Note that $\widetilde{X}(\lambda)$ and $X(\lambda)$ are necessarily different branches (see Figure~\ref{fig:schem-proof}) whereas $\widetilde{Y}(\mu)$ and $Y(\mu)$ may be the same, or distinct, branches (for example if $M$ overshoots a tipping point and returns to the same values after the upstream tipping occurs). Figure~\ref{fig:schem-proof} gives a schematic diagram that explains the hypothesis {\bf H3}. Note that in case of a tipping of the downstream system, the $\tilde{y}(\tilde{t})$ is a frozen tipping trajectory of the downstream system, whereas in case of tracking it is just a constant solution.

Finally, we can formulate a theorem that indicates the presence or absence of {\em Downstream within Upstream B-tipping} (DwUB-tipping in Figure~\ref{fig:tipping_seq}). 

\begin{theorem}[{\bf Condition for Downstream within Upstream B-tipping}]
\label{thm:criteria}
Suppose that {\bf H3} holds. Then there is an $r_*>0$ such that if $0<r<r_*$ then the upstream system will undergo B-tipping. In the following mutually exclusive cases, $0 < r_0 \leq r_*$ can be chosen such that:
\begin{itemize}
\item[1.]  If $M(\tilde{x}(t))\in \cL_{Y}$ for all $t$, then there is an $\epsilon_0(r)>0$ such that for any $0<r<r_0$ and $0<\epsilon<\epsilon_0(r)$ there will be downstream tracking for (\ref{eq:slowfastode2}).
\item[2.] If there are values $t_1 < t_2$ such that $M(\tilde{x}(t)) \notin\cL_Y$ for $t_1 < t < t_2$, then there is an $\tilde{\epsilon}_0(r)>0$ such that if $0<r<r_0$ and $0<\epsilon<\tilde{\epsilon}_0(r)$ then there is downstream B-tipping for (\ref{eq:slowfastode2}) within upstream tipping.
\end{itemize}
\end{theorem}

\proof
The conclusion for upstream B-tipping follows from parts 2 and 3 of Lemma~\ref{lem:Btip-up}, applying assumptions {\bf H3}. Statement 1 on downstream behaviour follows from noting there is an $\eta>0$ such that for $r$ small enough, {\bf H3} implies that $x^{[r]}(t)$ will either be within an $\eta$-neighbourhood of $X(\Lambda(s))$ for $s<s_f$, an $\eta$-neighbourhood of (a time-translation of) $\tilde{x}(t)$, or an $\eta$-neighbourhood of  $\tilde{X}(\Lambda(s))$ for $s>s_f$; see Figure~\ref{fig:schem-proof}. If all of these lie within $\cL_{Y}$, then tracking of the downstream system will follow for small enough $\epsilon$, where in general the largest $\epsilon_0$ depends on $r$.
Statement 2 on downstream B-tipping follows from noting that for all small enough $r$ to track the crossing of a threshold $\partial \cL_{Y}$, trajectories for small enough $\epsilon$ must switch to (and remain on) a new attractor %\todo{We don't have $Z(\lambda)$ any more so I've just removed it and left the new attractor unspecified...} 
by shadowing the frozen tipping trajectory $\tilde{y}(\tilde{t})$. 
\qed

\begin{corollary}
The requirement for sufficiently small $\epsilon$ in Theorem~\ref{thm:criteria} part 2 highlights that DwUB-tipping can only arise in the setting of an accelerating cascade.
\end{corollary}
In the following sections, we explore this and several other scenarios of tipping that can appear even in these simple cascades of tipping elements.

%\todo[inline]{This could be an alternative location for section 1.2?}

\section{Tipping scenarios in a cascade of two systems with hysteresis}
\label{sec:example}

To illustrate an application of the results in the previous section, we consider a class of example systems where both the upstream and downstream systems are simple one-dimensional systems (i.e., $n=1$, $m=1$)  We assume that both systems have hysteresis, but the coupling between the systems may take various forms. Specifically, we let upstream $x\in\R$ and downstream $y\in\R$ evolve according to
\begin{equation}
    \begin{aligned}
        \frac{\mathrm{d}x}{\mathrm{d}t}=& f(x,\Lambda(rt)),\\
        \epsilon\frac{\mathrm{d}y}{\mathrm{d}t}=&f(y,M(x)),
    \end{aligned}
    \label{eq:model}
\end{equation}
where $f(x,\lambda) : = 3x-x^3+\lambda$ has a hysteresis loop for $\lambda$ between $\lambda_l:=-2$ and $\lambda_u:=2$ and we set $r>0$ and $\epsilon>0$. The associated frozen system is therefore
\begin{equation}
	\begin{aligned}
		\frac{\mathrm{d}x}{\mathrm{d}t}=& f(x,\lambda),\\
		\epsilon\frac{\mathrm{d}y}{\mathrm{d}t}=&f(y,M(x)),
	\end{aligned}
	\label{eq:modelfrozen}
\end{equation}
There are stable branches of equilibria of the frozen upstream system at $x=X_{l}(\lambda)$ for $\lambda<\lambda_u$ and $x=X_{u}(\lambda)$ for $\lambda>\lambda_l$. These are connected by an unstable branch $x=U(\lambda)$ for $\lambda_l<\lambda<\lambda_u$, where $X_l(\lambda)<0<X_u(\lambda)$. Note that there are generic fold bifurcations at $\lambda_l$ and $\lambda_u$; hence, the system is bistable for $\lambda_l<\lambda<\lambda_u$.

We force the system (\ref{eq:model}) via the parameter shift
\begin{equation}
\begin{aligned}
\Lambda(s) &: = \lambda_-+[\lambda_+-\lambda_-]\frac{\tanh(s)+1}{2}
\label{eq:external_forcing}
\end{aligned}
\end{equation}
for real constants $\lambda_-$, $\lambda_+$, such that $\Lambda(s) \rightarrow \lambda_{\pm}$ as $s\rightarrow \pm\infty$ and $\Lambda(rt)$ is monotonic in $t$. We consider $r>0$ and the case  $\lambda_-<\lambda_u = 2$. 

\subsection{Upstream tipping behaviour}
\label{sec:upstream}

We first consider the behaviour of the upstream system. Here, for any $r>0$ there is a unique (pullback attracting) trajectory $x^{[r]}(t)$ limiting to the lower branch $X_l(\lambda_-)$ in the limit $t\rightarrow-\infty$, i.e. such that 
\begin{equation}
\lim_{t\rightarrow - \infty}x^{[r]}(t)=X_l(\lambda_-).
\label{eq:X-}
\end{equation}
We define (consistent with notation in section~\ref{sec:timescales_cascades})
$$
X_{-}:=\lim_{t\rightarrow -\infty}x^{[r]}(t) \mbox{ and }
X_{+}^{[r]}:=\lim_{t\rightarrow \infty}x^{[r]}(t).
$$
Using monotonicity of $\Lambda$, the form of $f$ and Lemma~\ref{lem:Btip-up}, we can classify the behaviour of $x^{[r]}(t)$ as follows:
\begin{itemize}
    \item {\bf Upstream tracking} If $\lambda_+<\lambda_u$ then for all $r>0$ there will be end-point tracking: $X_+^{[r]}=X_l(\lambda_+)$.
    \item {\bf Upstream B-tipping} If $\lambda_+>\lambda_u$ then for all $r>0$ there will be B-tipping: $X_+^{[r]}=X_u(\lambda_+)$.
\end{itemize}
Note that because $X_u$ and $X_l$ are basin stable (i.e. $X_{u,l}(\lambda)$ is in the basin of $X_{u,l}(\lambda')$ for any choices of $\lambda$ and $\lambda'$) results in \cite[Theorem 3.2]{ashwin2017parameter} about the lack of rate-induced tipping in the more restricted context of forward basin stability suggest rate-induced tipping from either branch is not possible, though rate-dependent effects such as safe overshoots may still be possible. Thus, due to this independence of the outcome on $r$, we can write $X_+=X_+^{[r]}$. 

\begin{figure}[!ht]
    \centering
    \includegraphics[width=0.5\textwidth]{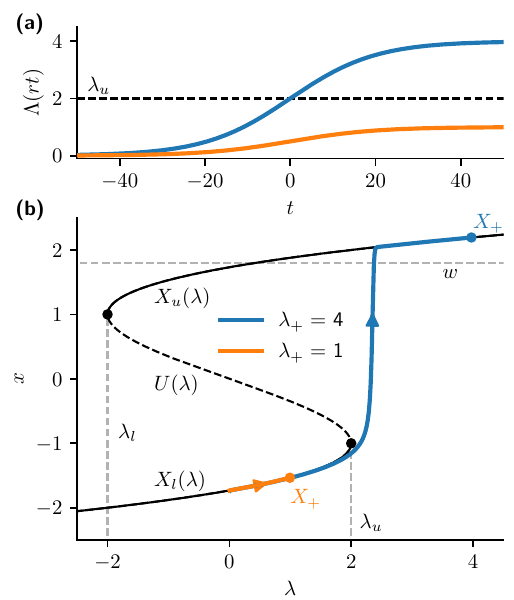}
    
    \caption{\textbf{Upstream system behaviour for ramp forcing} (a) Time series of external forcing profiles $\Lambda(rt)$ given by equation~\eqref{eq:external_forcing} with $\lambda_- = 0$ and $r = 0.05$. Blue forcing profile limits to $\lambda_+=4$, beyond the critical threshold, $\lambda_u=2$, for triggering tipping of the upstream system (horizontal black dashed line). Orange profile limits to $\lambda_+=1$, below the critical threshold. (b) Coloured curves display the response of upstream component $x$ of system~\eqref{eq:model} for the forcing profiles given in (a). Black solid/dashed curves correspond to the stable/unstable steady states, i.e. bifurcation diagram, for the frozen system~\eqref{eq:modelfrozen}. The blue trajectory crossing $\lambda_u$ signals the ``start time'' of tipping and crossing a threshold $w$ to enter a neighbourhood of the alternative state; we use this crossing to define the ``end time'' of tipping.}
    \label{fig:upstream_cases}
\end{figure}

Figure~\ref{fig:upstream_cases} provides examples for both cases listed above using the upstream component of \eqref{eq:model} with a slowly changing ($r=0.05$) external forcing, given by \eqref{eq:external_forcing}. The external forcing profiles, plotted in Figure~\ref{fig:upstream_cases}(a), limit to different levels $\lambda_+$, either above (blue), or below (orange) the critical threshold, $\lambda_u=2$.

The response of the upstream system, $x$, to this forcing is overlaid in colour on the bifurcation diagram of the frozen ($\Lambda(rt) = \lambda$) system~\eqref{eq:model} in Figure~\ref{fig:upstream_cases}(b). The orange trajectory demonstrates upstream tracking as the external forcing remains below $\lambda_u$, whereas, there is upstream B-tipping for the blue trajectory that exceeds the threshold $\lambda_u$, signalling the start of tipping, and crossing $w$ (the neighbourhood threshold of the alternative state) the end of tipping.

\subsection{Downstream behaviour}
\label{sec:downstream}

Following a similar approach, we can define different types of downstream behaviour. For this, we define a quantity that is the effective forcing of the downstream $y$ dynamics:
$$
\mu^{[r]}(t):=M(x^{[r]}(t)),
$$
and define 
$$
M_-:=\lim_{t\rightarrow-\infty} \mu^{[r]}(t).
$$
This is independent of $r$ in similarity to assumption (\ref{eq:X-}) about $X_-$.
We assume that $M_-<\lambda_u$
and so for any $r>0$, $\epsilon>0$ there is a unique (pullback attracting) trajectory $y^{[r,\epsilon]}(t)$ that limits to $X_l(M_-)$. Specifically, the limit
$$
Y_-:=\lim_{t\rightarrow -\infty} y^{[r,\epsilon]}(t).
$$
satisfies $Y_-=X_l(M_-)$. We define
$$
Y_+^{[r,\epsilon]}:=\lim_{t\rightarrow +\infty} y^{[r,\epsilon]}(t).
$$
and
$$
M_{+}:=\lim_{t\rightarrow +\infty} \mu^{[r]}(t)
$$
(where we note that no $r$ dependence is needed as $X_+ = X_+^{[r]}$).

Now we are in a position to classify the behaviour of the downstream system $y^{[r,\epsilon]}(t)$ as follows:
\begin{itemize}
    \item {\bf Downstream tracking:} If $\mu^{[r]}(t)<\lambda_u$ for all $t>0$ then we have end-point tracking: $Y_+^{[r,\epsilon]}=X_l(M_+)$.
    \item {\bf Downstream B-tipping:} If $M_+>\lambda_u$ then we have B-tipping: $Y_+^{[r,\epsilon]}=X_u(M_+)$.
    \item {\bf Downstream overshoot
    :} If $M_+<\lambda_u$ and there is an $r$ and a $t$ (depending on $r$) such that $\mu^{[r]}(t)>\lambda_u$, then for small enough $\epsilon$ there is B-tipping for the downstream system. For large enough $\epsilon$ there is end-point tracking.
\end{itemize}
Observe that downstream overshoot can only appear if $\mu^{[r]}(t)$ is non-monotonic, where $M(x^{[r]}(t))$ briefly overshoots the threshold $\lambda_u$ and then returns. For a sufficiently large timescale separation $\epsilon$, one can have an overshoot without causing tipping \cite{ritchie2019inverse}. The timings of when the upstream and downstream tipping occur can be classified according to Figure~\ref{fig:tipping_seq}.

\subsection{The dependence of tipping scenario on coupling}
\label{sec:depend}

Clearly, the coupling $M$ can be of a variety of forms, and some of the methods above can be applied independently of the form of $M$. However, for definiteness, we consider in the following sections two specific cases of coupling $M(x)=M_i(x)$ for $i=1,2$  defined as in Table~\ref{tab:cases}. We explore the downstream tipping scenarios for (\ref{eq:model}). We make the standing assumptions
\begin{equation}
\lambda_-<\lambda_u=2, ~~X_-=X_l(\lambda_-),~~Y_-=X_l(M_-).
\label{eq:standingassumptions}
\end{equation}
and so $M_-<\lambda_u$.
Moreover, from hereon we will assume $\lambda_+>\lambda_u$, such that the upstream system undergoes B-tipping (blue trajectory in Figure~\ref{fig:upstream_cases}) and so $x^{[r]}(t)$ limits to $X_+=X_u(\lambda_+)$ in forward time. Specifically, we choose the default parameters as in Table~\ref{tab:cases}
that determine the forcing of the upstream system. 

The classification of a scenario according to timing, as in Figure~\ref{fig:tipping_seq}, is dependent on the thresholds used to define the onset and offset times of the tipping; respectively (for example  $\lambda=\lambda_u$, $x=w$ for upstream, and $M=\lambda_u$, $y=w$ for downstream tipping as in Figure~\ref{fig:upstream_cases}). We say that tipping has occurred only if the offset threshold is passed: the system can overshoot the onset threshold but not reach the offset threshold. Note that other definitions of threshold are possible, for example, onset and offset could be based on the speed of evolution exceeding some threshold, or more generally (as in \cite{ashwin2025early}) by setting thresholds for some observable that is large when tipping is underway. We discuss in Appendix~\ref{app:boundaries} the numerical methods used to estimate the boundaries between different tipping scenarios as well as present a sensitivity study on the offset threshold choice $w$.

\begin{table}[!ht]
\centering
\caption{\label{tab:cases} Cases for coupling function $M(x)$ that are considered for the system (\ref{eq:model}). In all cases, $a$, $b$, $c$, $d$ are constants, $X_-$ is a stable equilibrium of the past limit system, and we choose $x_0=X_-$ as the initial condition. We take other default values $r=0.05$, $\epsilon=0.05$, $\lambda_-=0$, $\lambda_+=4$, $x_0=X_-=-\sqrt{3}$ and $w=1.8$.}

~

\begin{tabular}{l|l|l}
Name & functional form & 
default values\\
\hline
Linear coupling & $M_1(x):=a+b(x-X_-)$ 
& $a=0$, $b=1$\\
Localised coupling & $M_2(x):=a+b\sech (c(x-d))$ 
& $a=0$, $b=1$, $c=2$, $d=0.5$\\
\hline
\end{tabular}
\end{table}

\section{Linear coupling}
\label{sec:linear}

For the case of linear coupling, $M(x)=M_1(x)=a+b(x-X_-)$, we note that 
$M_-=a$, $M_+=a+b(X_+-X_-)$
and so $M_-<\lambda_u$ implies
$a<\lambda_u$. For this scenario, noting that $M(x(t))$ will be monotonic for our choice of $f$, so for small enough $\epsilon$ we will have downstream tracking if $M_+<\lambda_u$
and downstream B-tipping if $M_+>\lambda_u$.

\subsection{Tipping behaviour for linear coupling}
\label{sec:lin_coupling}

\begin{figure}[!ht]
    \centering
    \includegraphics[width=\linewidth]{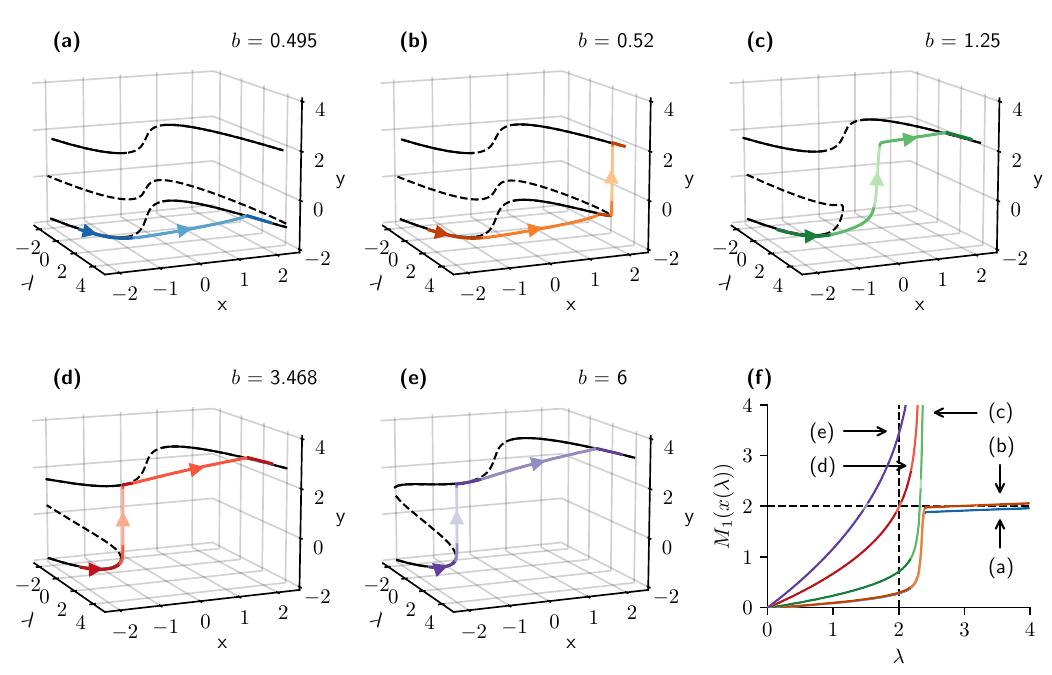}
    
    \caption{\textbf{Tipping behaviour for linear coupling} (a)--(e) 3D bifurcation diagrams for system~\eqref{eq:model} with linear coupling (black curves). System response (coloured curves) overlaid on top, for different coupling strengths $b$. (f) Plot showing order in which critical thresholds are crossed for the upstream system ($\lambda=2$) and downstream system ($M_1(x(\lambda)) = 2$) for different coupling strengths used in (a)--(e). Parameter values as given in Table~\ref{tab:cases}, except for $b$ which are provided for each 3D panel.}
    \label{fig:3d_bif_diag_lin}
\end{figure}

We choose default parameters as in Table~\ref{tab:cases} that determine the linear coupling from the upstream to the downstream system. We study how the coupling strength, $b$, promotes different behaviour. Here, we choose $\epsilon$ small, such that we consider the setting of an accelerating cascade of tipping elements. Five qualitatively different scenarios are shown in Figure~\ref{fig:3d_bif_diag_lin} (a)--(e) depending on the value of $b$.
In these figures, we show (black lines) bifurcation diagrams with stable and unstable branches for the combined frozen system (\ref{eq:modelfrozen}). We also show a trajectory (coloured line) representing a pullback attractor for the nonautonomous system (\ref{eq:model}) with the slowly varying parameter shift (\ref{eq:external_forcing}). The components limit, in the past, to the lower branch stable equilibrium for both subsystems. This is a dynamic trajectory, so we try to indicate the direction of time with an arrow and the regions of rapid change (associated with tipping) by using lighter colours. The bifurcation structures for both the upstream and downstream systems are the same, meaning that the upstream (resp. downstream) system cross a fold bifurcation when the forcing reaches $\lambda=2$ (resp.  $M_1(x(\lambda)) = 2$). Figure~\ref{fig:3d_bif_diag_lin}(f) plots the profiles of $(\Lambda(rt),M_1(x(\Lambda(rt))))$ for different coupling strengths $b$. The fold bifurcation at $\lambda=2$ (marked by the vertical black dashed line) is always crossed as the same external forcing profile is used -- blue forcing profile in Figure~\ref{fig:upstream_cases}(a). However, if and when the fold bifurcation at $M_1(x(\lambda)) = 2$ is crossed, changes based on the coupling strength $b$, as illustrated in Figure~\ref{fig:3d_bif_diag_lin}(f).

If the coupling is sufficiently weak, then the downstream system never crosses the fold bifurcation at $M_1(x(\lambda)) = 2$, as illustrated by the blue curve in Figure~\ref{fig:3d_bif_diag_lin}(f) and verified by applying Theorem~\ref{thm:criteria}, case 1. Consequently, only the upstream system tips while the downstream system tracks, corresponding to UB-tipping. We can visualise this scenario by superimposing the system trajectory (coloured blue) on the bifurcation diagram for the full system~\eqref{eq:model} in Figure~\ref{fig:3d_bif_diag_lin}(a). 

The system is initialised on the lower left (as viewed) stable branch. Upon increasing the forcing, $\Lambda(rt)$, the system tracks this lower branch until a fold bifurcation is reached. Crossing the fold bifurcation causes the upstream system to tip and travel faster (indicated by the lighter shading) until the alternative lower stable branch on the right is reached (corresponding to a tipped upstream system but a tracked downstream system). This branch also features a fold bifurcation (associated with tipping of the downstream system). Due to the continued increase in forcing, the system moves towards this second fold bifurcation. However, the forcing stabilises before this second fold bifurcation is reached and therefore does not undergo tipping again.

If the forcing had continued, or the coupling strength was a little stronger, downstream B-tipping would occur after the upstream B-tipping, corresponding to DaUB-tipping. The orange curves in Figure~\ref{fig:3d_bif_diag_lin}(b) and (f) illustrate this scenario by crossing the second bifurcation, caused by a small increase in coupling strength. Notice that after crossing the second fold bifurcation, the downstream B-tipping occurs faster (denoted by the lightest shading of orange) than the previous upstream B-tipping (intermediate shading) due to this being an accelerating cascade. If instead the downstream threshold is crossed while the upstream system is undergoing B-tipping (i.e. for stronger coupling), then DwUB-tipping becomes possible, as shown by the green trajectories in Figure~\ref{fig:3d_bif_diag_lin}(c) and (f). The fold bifurcation at $\lambda = 2$ is still crossed first, causing upstream B-tipping to commence first. However, as seen in Figure~\ref{fig:3d_bif_diag_lin}(c), downstream B-tipping then begins and subsequently finishes (green shading temporarily changes to lightest shading) before upstream B-tipping finishes.

There is a qualitative change in behaviour if the coupling strength is sufficiently strong such that the thresholds for the upstream and downstream systems are crossed simultaneously; see the red trajectory in Figure~\ref{fig:3d_bif_diag_lin}(f). Once again, due to the much faster timescale on the downstream system, the downstream system tips rapidly causing the system to essentially land on another fold bifurcation, Figure~\ref{fig:3d_bif_diag_lin}(d), after which tipping of the upstream system occurs at the slower speed (intermediate red shading). The system then again converges to the stable branch on the upper right. 

For very strong coupling, the system may cross the fold bifurcation corresponding to the downstream system ($M_1(x(\lambda)) = 2$) before the fold bifurcation ($\lambda = 2$) of the upstream system is exceeded, see purple curve in Figure~\ref{fig:3d_bif_diag_lin}(f). This leads to UaDB-tipping (upstream B-tipping after downstream B-tipping), see Figure~\ref{fig:3d_bif_diag_lin}(e). Here, the system crosses its first fold bifurcation, and downstream B-tipping occurs to an intermediate stable branch on the upper left. As the forcing, $\lambda$, continues to increase, the second fold is crossed, and upstream B-tipping at a slower speed (middle shading). The sequence of events in this scenario can give the misleading impression that the downstream B-tipping causes the upstream B-tipping, but clearly this cannot be the case as there is no coupling from the downstream to the upstream system in \eqref{eq:model}. 

\begin{figure}[!ht]
    \centering
    \includegraphics[width=\textwidth]{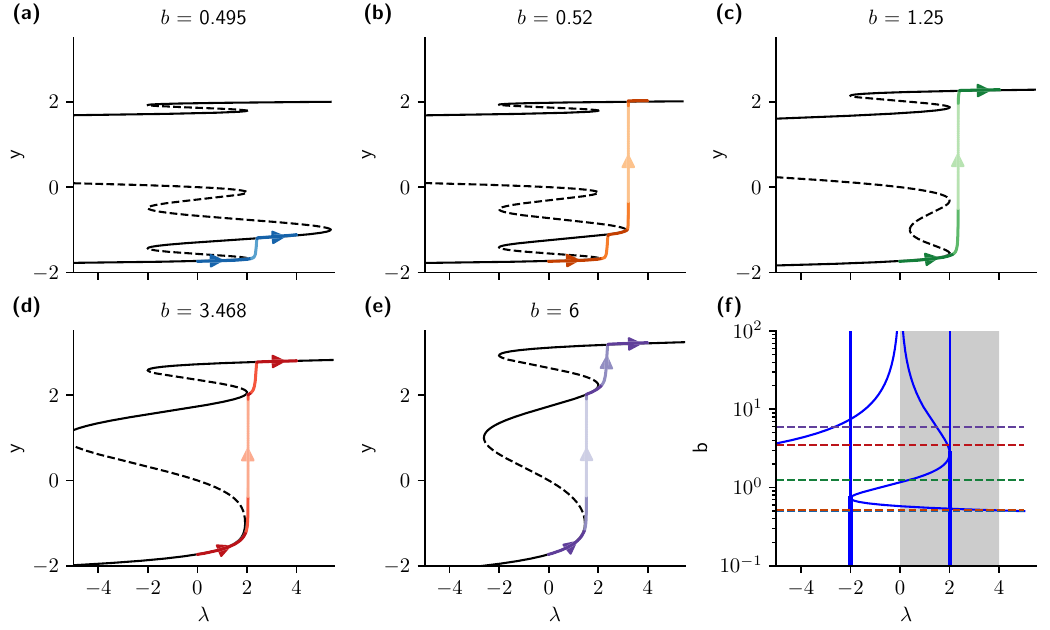}
    
    \caption{\textbf{Bifurcation analysis for linear coupling} (a)--(e) Bifurcation diagrams plotted in the $(\lambda,y)$ projection for system~\eqref{eq:model} with linear coupling $M_1(x)$ (black curves). System response (coloured curves) overlaid on top, for the different coupling strengths, $b$. (f) Two-parameter bifurcation diagram in the $(\lambda,b)$ plane showing the location of fold bifurcations (blue). Thick blue lines correspond to three fold bifurcations at $\lambda=-2$ and $\lambda=2$ for sufficiently weak coupling strengths. The grey shaded region indicates the parameter shift range of $\lambda$. Horizontal dashed coloured lines indicate coupling strengths, $b$, used in (a)--(e). Parameter values as given in Table~\ref{tab:cases}, except for $b$, which are provided for each one-parameter bifurcation diagram.}
    \label{fig:2par_bif_diag_lin}
\end{figure}

Figure~\ref{fig:2par_bif_diag_lin}(a)--(e) provides further clarity on how the sequences of bifurcations in the $(\lambda,y)$-plane change for different coupling strengths $b$. Recall that the stability of each equilibrium is determined by the full 2D Jacobian of \eqref{eq:slowfastodefrozen}, which can lead to counter-intuitive structures in this projection, such as apparent adjacent unstable branches. However, these equilibria might be located at very different $x$ values. %\todo{Check! Added this for Ref 1, com 2} 
Figure~\ref{fig:2par_bif_diag_lin}(f) provides the 2-parameter bifurcation diagram plotting the location of the fold bifurcations in the $(\lambda, b)$-plane.

For weak coupling strengths, there are up to nine steady states, separated by up to eight fold bifurcations. In Figure~\ref{fig:2par_bif_diag_lin}(a), and the weakest coupling strengths, a fold bifurcation is crossed at $\lambda = 2$, causing a small jump in this projection due to upstream B-tipping (as seen in Figure~\ref{fig:3d_bif_diag_lin}(a)). The system subsequently moves towards the fold bifurcation furthest right, this is associated with additional downstream B-tipping. However, this fold is beyond $\lambda=4$, and therefore the system stabilises before the fold is reached. For a higher coupling strength, $b$, this fold moves left (compare Figures~\ref{fig:2par_bif_diag_lin}(a) and (b)). This is represented in Figure~\ref{fig:2par_bif_diag_lin}(f) by the blue curve that comes in from the right-hand side. The grey shaded region represents the range over which the parameter $\lambda$ is varied. Hence, the intersection of the blue curve with the right edge of the grey shaded region ($\lambda = \lambda_+ = 4$) provides the minimum coupling strength ($b_{\min}\approx0.51$) required for a tipping cascade.

For sufficiently strong coupling, such that the fold moves below $\lambda = 2$, no steady states exist beyond $\lambda = 2$, except for the uppermost state, corresponding to the tipped upstream and downstream state. Therefore, this coupling strength approximately (given the choice of threshold $w$) marks the beginning of downstream within upstream B-tipping in the limits of $r,\epsilon\downarrow 0$.

The fold continues to move backwards towards $\lambda = -2$ for increasing $b$, while simultaneously the neighbouring folds, either side, move closer together in the $y$ direction. Once the fold reaches $\lambda = -2$, all three folds collide to leave only one fold bifurcation: this corresponds to a degenerate higher-codimension bifurcation where the degeneracy is due to the one-directional coupling. As the coupling strength, $b$, continues to increase, this fold moves to higher values of $\lambda$. However, this coalescence has no impact on the tipping behaviour, which remains downstream within upstream B-tipping, as illustrated by the green curve in Figure~\ref{fig:2par_bif_diag_lin}(c). The system crosses the fold bifurcation at $\lambda = 2$ and tips straight to the uppermost stable state without reaching an intermediate steady state.

The same dynamics occur as before for increasing the coupling strength, $b$, further. The fold moves to higher values of $\lambda$ and the two neighbouring folds, either side, located at $\lambda = 2$, move closer together in the $y$ direction. Once the newly formed fold reaches $\lambda = 2$, the three folds collide and again only one fold remains. Notice that the special case (Figure~\ref{fig:2par_bif_diag_lin}(d)) of passing the critical thresholds for the upstream and downstream systems simultaneously does not occur precisely at this collision; see the red dashed line in Figure~\ref{fig:2par_bif_diag_lin}(f). Instead, due to transient effects, the location of this upper boundary for downstream within upstream B-tipping will be for a slightly greater coupling strength, the degree of which depends on the timescale separation, $\epsilon$. Similar arguments determine that the lower boundary for downstream within upstream B-tipping does not coincide exactly with the intersection of the fold curve emanating from the right and the three folds line at $\lambda=2$ in Figure~\ref{fig:2par_bif_diag_lin}(f), but a little below (though still above the dashed orange line). 

Increasing the coupling strength further pushes the fold to lower values of $\lambda$ and asymptotes to $\lambda = 0$ as $b\rightarrow\infty$. For these strong coupling strengths, the movement of this fold further below $\lambda = 2$ enables the possibility of the downstream system to tip (and reach an intermediate stable state) before the upstream system tips at $\lambda = 2$, see Figure~\ref{fig:2par_bif_diag_lin}(e). Although not shown, if $\lambda_+$ is chosen such that the system crosses the first fold but does not reach $\lambda = 2$ then only the downstream system tips. Another fold curve in Figure~\ref{fig:2par_bif_diag_lin}(f) comes in from the left (that is associated with downstream B-tipping from the upper branch to the lower branch) and also asymptotes to $\lambda = 0$. However, this has little impact on the upstream tipping behaviour. 

\subsection{Tipping regimes for linear coupling}

\begin{figure}[!ht]
    \centering
    \includegraphics[width=\textwidth]{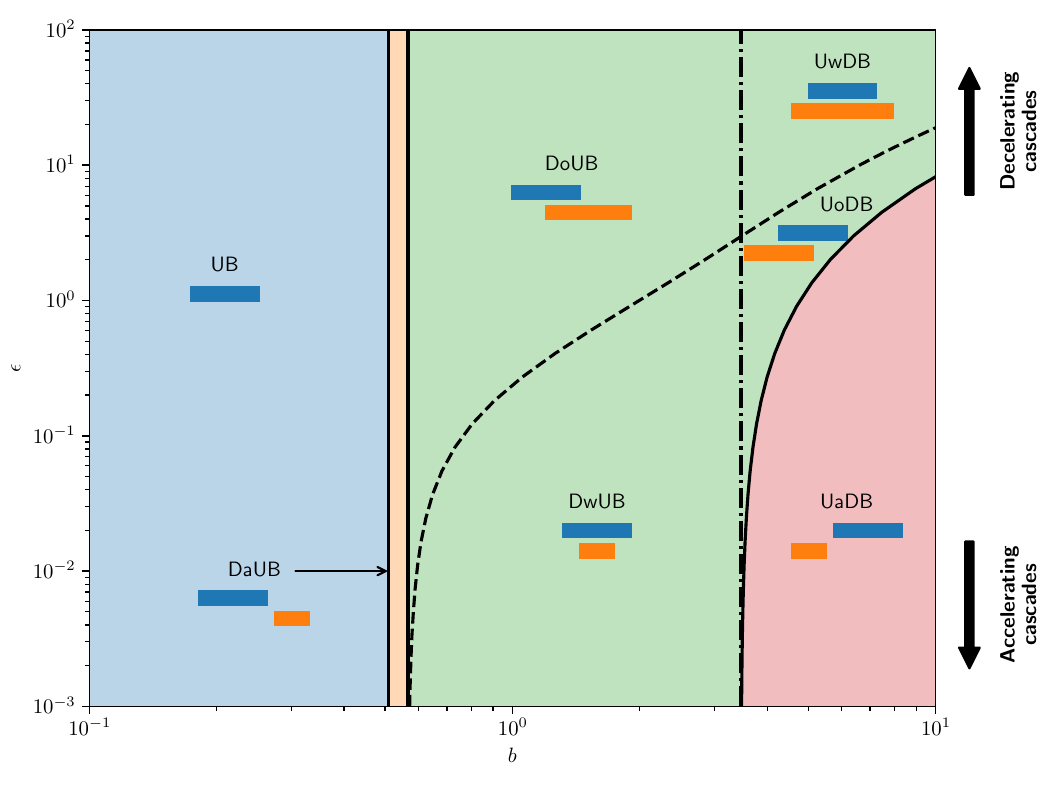}
    
    \caption{\textbf{Tipping regimes for linear coupling} Tipping regimes for system~\eqref{eq:model} with linear coupling $M_1(x)$ dependent on the coupling strength $b$, and the timescale separation, $\epsilon$. Blue region (UB): upstream B-tipping only (downstream tracking). Orange region (DaUB): downstream B-tipping after upstream B-tipping. The green region represents when there is at least a brief period of both upstream and downstream systems tipping. This region is partitioned by the relative timings of the tipping onset (vertical dash-dot line) and offset (curved dashed line) for the upstream and downstream systems. This gives the regions (DwUB): downstream within upstream B-tipping; (DoUB): downstream overlap upstream B-tipping; (UoDB): upstream overlap downstream B-tipping; (UwDB): upstream within downstream B-tipping. Red region (UaDB): upstream after downstream B-tipping that reaches an intermediate state. The boundaries are computed as outlined in Section~\ref{sec:depend}. The figure labels use the abbreviations UB, DaUB, DwUB, DoUB, UoDB, UwDB, and UaDB; throughout the text we refer to the corresponding tipping regimes as UB-tipping, DaUB-tipping, etc. Coloured blocks illustrate the tipping duration of the upstream (blue) and downstream (orange) systems as in Figure~\ref{fig:tipping_seq}. Parameter values as given in Table~\ref{tab:cases}.}
    \label{fig:tipping_regimes_lin}
\end{figure}

We summarise, in terms of the classification in Figure~\ref{fig:tipping_seq}, the different tipping regimes that can be found in system~\eqref{eq:model} with linear coupling in Figure~\ref{fig:tipping_regimes_lin} for different values of the timescale separation $\epsilon$ and coupling strength $b$. A horizontal cross-section in the lower half of the figure corresponds to an accelerating cascade of tipping elements as considered up to now. Moving from left to right along such a cross-section, weak coupling first gives rise to UB-tipping (labelled UB; blue region), in which only the upstream system undergoes B-tipping. Over a narrow interval of coupling strengths, DaUB-tipping (labelled DaUB; orange region) occurs, where downstream B-tipping occurs after upstream B-tipping. This is followed by a region of predominantly DwUB-tipping (labelled DwUB; green region), in which downstream B-tipping occurs within the period of upstream B-tipping. The vertical dash-dot line separates coupling strengths at which the order of tipping onset (defined as crossing the forcing threshold) reverses between the upstream and downstream systems. For sufficiently strong coupling strengths, the downstream threshold is crossed first, and, due to the faster dynamics, the downstream system tips before the upstream threshold is crossed. Therefore, UaDB-tipping (labelled UaDB; red region) occurs, in which upstream B-tipping happens after downstream B-tipping and a different intermediate state (compared with the DaUB-tipping regime) is temporarily reached. The boundaries are computed as outlined in Section~\ref{sec:depend}.

The coloured bars used to indicate the duration of tipping for the upstream (blue) and downstream (orange) systems graphically illustrate these different regimes. The tipping dynamics for the upstream system is independent of $b$ and $\epsilon$ and therefore the tipping onset and duration (width of the blue box) do not change. For increasing $b$, but fixed $\epsilon$, the onset of downstream tipping moves earlier, but the duration remains largely unchanged. Therefore, the orange bar starts to the right of the blue bar, corresponding to DaUB-tipping (labelled DaUB), then shifts left, briefly overlapping with the blue bar (DoUB-tipping), before becoming completely contained within the blue bar (DwUB-tipping). Beyond the vertical dashed line the orange bar starts first and so there is a brief region of overlap the other way around (UoDB-tipping), and then for strong coupling the orange bar stops before the blue bar starts (UaDB-tipping).          

The regions change slightly if considering a decelerating cascade of tipping elements (horizontal cross section near top of Figure~\ref{fig:tipping_regimes_lin}). Although, the boundaries marking the end of UB-tipping and DaUB-tipping (two vertical black solid lines close together), and the boundary for the simultaneous onset of upstream and downstream tipping (vertical dash-dot line) are independent of the timescale separation. 
%However, the DaUB region changes from orange (small $\epsilon$) to green (large $\epsilon$), indicating that an intermediate state is no longer reached for decelerating cascades. Since the dynamics of the downstream system is much slower, the effective forcing on the downstream system during the upstream B-tipping is now fast. Therefore, despite not crossing the downstream threshold during the upstream B-tipping, the downstream system is forced out of equilibrium. 
%The downstream system then does not have time to converge back to the equilibrium branch before it disappears as the threshold is subsequently crossed, meaning that an intermediate state is never reached. 
%Note that increasing $\epsilon$ does not change the onset timing of downstream tipping and hence the dash-dot line indicating the simultaneous onset of tipping is vertical. 
However, the simultaneous offset boundary does vary based on the timescale separation as indicated by the curved black dashed line.
%indicates when the upstream and downstream systems finish tipping at the same time. 
Hence, for a fixed $b$ corresponding to DwUB-tipping at small $\epsilon$, 
%given onset of downstream tipping (determined by $b$) 
sufficiently large $\epsilon$ eliminates DwUB-tipping, with the tipping of the two components instead overlapping (DoUB-tipping). Additionally, for sufficiently large $\epsilon$, the duration of the downstream tipping exceeds that of the upstream tipping. Thus, for coupling strengths beyond the vertical black dash-dot line, UwDB-tipping becomes possible, with upstream B-tipping completely contained within downstream B-tipping.    

\section{Localised coupling}
\label{sec:localised}

In this section, we consider the case $M(x)=M_2(x)=a+b\sech (c(x-d))$ corresponding to a coupling that is localised around $x=d$, where we choose $d$ to be a transient point during the tipping of the upstream system. If $c$ is sufficiently large then $M_{\pm} \approx a$ while $M(d) = a + b$. That is, there is a turning point at $x = \tilde{x}:= d$. We can classify the behaviour of the downstream system in terms of $M_{\pm}$ and $M(\tilde{x})$.

% For more general types of continuous coupling $M(x)$ if there is a global maximum $\tilde{x}$ on the interval $(X_-,X_+)$ then the downstream behaviour may relate to the upstream as in the linear coupling case but in addition there may be overshoots of thresholds in the downstream system.

\subsection{Tipping behaviour for localised coupling}

\begin{figure}[!ht]
    \centering
    \includegraphics[width=\linewidth]{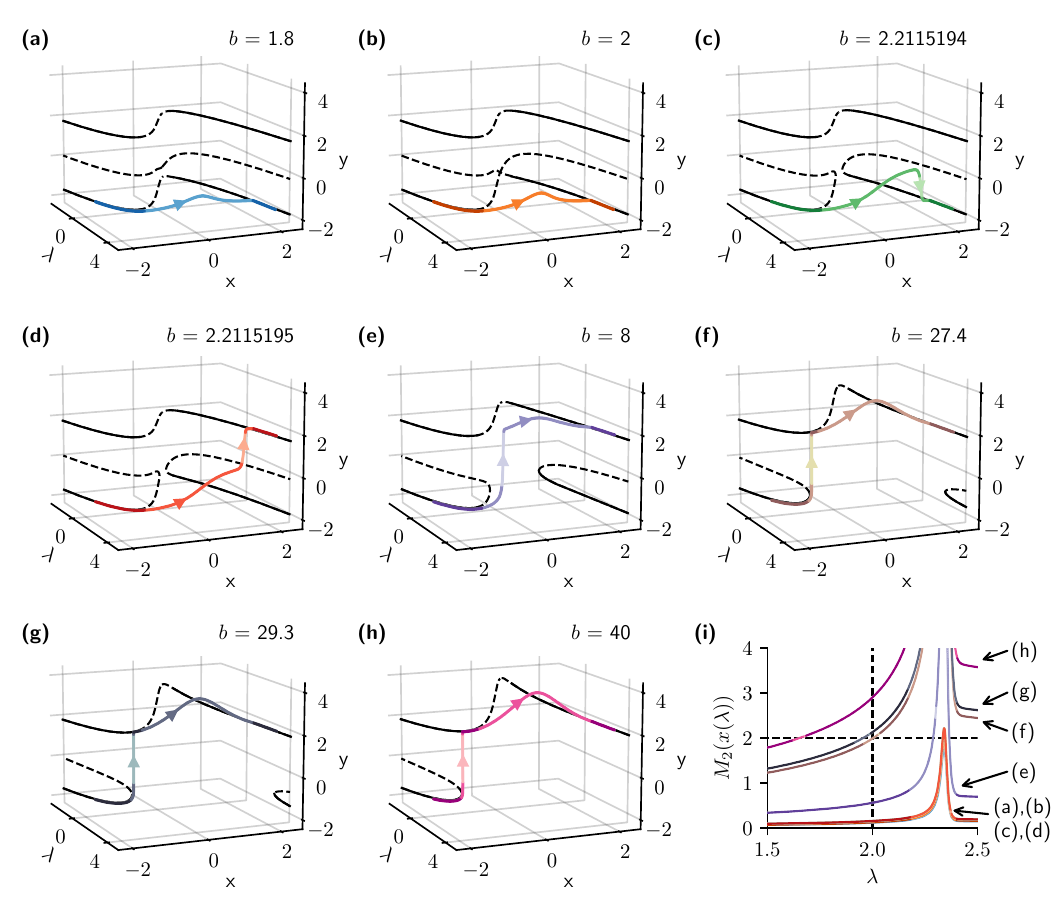}
    
    \caption{\textbf{Tipping behaviour for localised coupling} (a)--(h) 3D bifurcation diagrams for system~\eqref{eq:model} with localised coupling (black curves). System response (coloured curves) overlaid on top, for different coupling strengths $b$. (i) Plot showing order in which critical thresholds are crossed for the upstream system ($\lambda(t)=2$) and downstream system ($M_2(x(\lambda(t))) = 2$) for different coupling strengths used in (a)--(h). Parameter values as given in Table~\ref{tab:cases}, except for $b$, which are provided for each 3D panel.}
    \label{fig:3d_bif_diag_loc}
\end{figure}

We consider the response of the following system for localised coupling. For fixed $a$, $b$, $c$, and $d$ there may or may not be an excursion of $M(x)$ that can take $y$ over the downstream tipping threshold, i.e., depending on whether $M(\tilde{x})$ is above or below $\lambda_u$. We choose the default parameters as in Table~\ref{tab:cases} that determine the coupling from upstream to downstream systems, and vary the coupling strength $b$ (other settings are standard as per \eqref{eq:standingassumptions}). Thus $M_\pm = 0$ and $M(\tilde{x}) = b$. 
In Figure~\ref{fig:3d_bif_diag_loc} we report the qualitatively different behaviours of system~\ref{eq:model} for different coupling strengths, $b$, of the coupling function.

If the coupling strength $b$ is sufficiently small (blue trajectory), then the threshold of the downstream system is never reached. This can be seen in the 3D bifurcation diagram in Figure~\ref{fig:3d_bif_diag_loc}(a), where only upstream B-tipping occurs. Applying Theorem~\ref{thm:criteria}, case 1, confirms that the downstream system remains tracking. This scenario of UB-tipping, with upstream B-tipping and downstream tracking, is also observed in Figure~\ref{fig:3d_bif_diag_loc}(b), where the threshold for the downstream is only touched before the forcing moves away. 

For stronger coupling $b$, a brief overshoot of the downstream threshold occurs such that there are now two crossings of the downstream threshold, seen by plotting the profiles of $(\Lambda(rt),M_2(x(\Lambda(rt))))$ in Figure~\ref{fig:3d_bif_diag_loc}(i). 
Figure~\ref{fig:3d_bif_diag_loc}(c) illustrates that once the fold bifurcation at $\lambda=2$ is crossed, the system leaves the base state. The green trajectory starts to tip in the upstream $x$ direction, while also increasing in the downstream $y$ direction. However, the effective forcing on the downstream system is reversed sufficiently quickly (for the distance of overshoot) such that the trajectory does not tip to the upper stable branch. Therefore, this remains UB-tipping, with downstream tracking. For yet stronger coupling $b$, the red trajectory deviates too far and consequently the downstream system tips to the upper stable branch; see Figure~\ref{fig:3d_bif_diag_loc}(d). At even stronger coupling strengths, the two lower stable branches move further apart, resulting in DwUB-tipping, where downstream B-tipping occurs within upstream B-tipping, as shown in Figure~\ref{fig:3d_bif_diag_loc}(e).

The brown trajectory in Figure~\ref{fig:3d_bif_diag_loc}(i) displays the special case of crossing the thresholds for the upstream and downstream systems simultaneously. In the 3D bifurcation (Figure~\ref{fig:3d_bif_diag_loc}(f)), the lower fold is crossed and again, because of the timescale separation, downstream B-tipping occurs, and the trajectory lands close to the upper fold, leading to upstream B-tipping. For larger coupling strengths, the downstream threshold will be crossed first, see grey and pink trajectories in Figure~\ref{fig:3d_bif_diag_loc}(i). There is a range of $b$ corresponding to an overshoot where there is bistability when the forcing is stabilised (Figure~\ref{fig:3d_bif_diag_loc}(g)). However, for the largest coupling strengths, the bistability is lost, meaning that an overshoot no longer occurs (Figure~\ref{fig:3d_bif_diag_loc}(h)). Additionally, crossing the downstream threshold sufficiently long before the upstream threshold (pink trajectory in Figure~\ref{fig:3d_bif_diag_loc}(i)) enables the downstream system to B-tip to an intermediate stable state in the top left of Figure~\ref{fig:3d_bif_diag_loc}(h), resulting in UaDB-tipping (upstream B-tipping after downstream B-tipping). 

\begin{figure}[!ht]
    \centering
    \includegraphics[width=\linewidth]{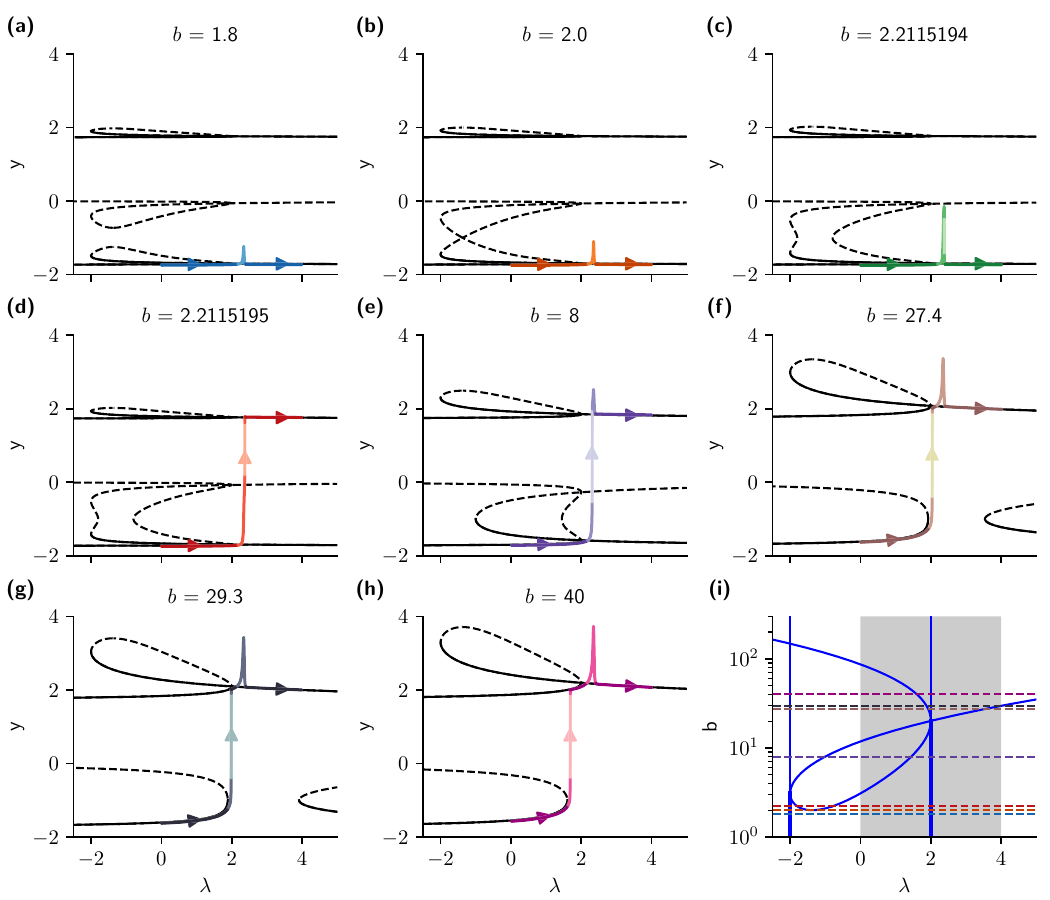}
    
    \caption{\textbf{Bifurcation analysis for localised coupling} (a)--(h) Bifurcation diagrams plotted in the $(\lambda,y)$ projection for system~\eqref{eq:model} with localised coupling $M_2(x)$ (black curves). System response (coloured curves) overlaid on top, for the different coupling strengths, $b$. (i) Two-parameter bifurcation diagram in the $(\lambda,b)$ plane showing the location of fold bifurcations (blue). Thick blue lines correspond to three fold bifurcations at $\lambda=-2$ and $\lambda=2$ for sufficiently weak coupling strengths. Grey shaded region indicates the parameter shift range of $\lambda$. Horizontal dashed coloured lines indicate coupling strengths $b$ used in (a)--(h). Parameter values as given in Table~\ref{tab:cases}, except for $b$, which are provided for each one-parameter bifurcation diagram.}
    \label{fig:2par_bif_diag_loc}
\end{figure}

Figure~\ref{fig:2par_bif_diag_loc}(a)--(h) plots the bifurcation diagram in the $(\lambda,y)$-plane for different levels of the coupling strength, $b$. The two-parameter bifurcation diagram, showing the location of the fold bifurcations, is given in Figure~\ref{fig:2par_bif_diag_loc}(i). For the weakest coupling strengths $b$ (e.g. blue dashed line in Figure~\ref{fig:2par_bif_diag_loc}(i)), there are six fold bifurcations (three at both $\lambda=\pm2$) giving rise to nine branches of steady states, as shown in Figure~\ref{fig:2par_bif_diag_loc}(a). The system starts on the lower stable branch, and as the forcing increases, the blue trajectory reaches a fold at $\lambda=2$. At the fold, the original stable state disappears, following a collision with the unstable branch, and this causes a small spike in the trajectory that corresponds to upstream B-tipping. However, another stable state is nearby (in the $y$ direction), and the system converges to this alternative stable state.  

For $b=2$, the two lower loops of steady states seen in Figure~\ref{fig:2par_bif_diag_loc}(a) touch, forming a cusp bifurcation, see Figure~\ref{fig:2par_bif_diag_loc}(b). In this projection, the orange trajectory is similar to the blue trajectory, only with a slightly larger spike. Note that the orange dashed line in Figure~\ref{fig:2par_bif_diag_loc}(i) is tangential to the fold bifurcation curve and therefore gives the cusp bifurcation. Thus in Figure~\ref{fig:2par_bif_diag_loc}(c), by increasing $b$, two new fold bifurcations are created, which move further apart in the $\lambda$ direction. The green trajectory now produces a large spike, but the downstream system still avoids transitioning to the upper stable branch as the coupling strength is just below the critical level. In Figure~\ref{fig:2par_bif_diag_loc}(d), the coupling strength $b$ is slightly increased to just above the critical level, which causes downstream B-tipping to the upper stable branch despite no qualitative change in the bifurcation structure. 

{\em Figure~\ref{fig:2par_bif_diag_loc} highlights the rate dependency on B-tipping, and specifically the critical coupling strength will change based on the timescale separation between the two systems, $\epsilon$, as well as the rate parameter, $r$, in the external forcing, $\Lambda(rt)$.}

In similarity to the case for linear coupling, the three folds that are close to each other near $\lambda=-2$ in Figure~\ref{fig:2par_bif_diag_loc}(d), move closer together upon increasing $b$. They collide, leaving only one fold, when the thick line at $\lambda=-2$ in Figure~\ref{fig:2par_bif_diag_loc}(i) changes to normal thickness and is tangent to the other fold curve. The remaining fold then moves in the positive $\lambda$ direction on increasing $b$. Simultaneously, the fold close to $\lambda=0$ in Figure~\ref{fig:2par_bif_diag_loc}(d) moves towards the two folds at $\lambda=2$, see Figure~\ref{fig:2par_bif_diag_loc}(e). Little change is observed in the purple trajectory, except for a small spike in the downstream response when reaching the alternative stable state before converging.

For stronger coupling $b$, the three folds close to $\lambda=2$ and to each other collide when the bold line at $\lambda=2$ in Figure~\ref{fig:2par_bif_diag_loc}(i) terminates and is tangential to the other fold curve. One fold remains, that retreats to smaller values of $\lambda$, while the other fold continues to move to larger values of $\lambda$. Figure~\ref{fig:2par_bif_diag_loc}(f) corresponds to the special case of crossing the folds for the upstream and downstream systems simultaneously. The first fold is crossed, causes downstream B-tipping due to the fast timescale, and then the solution lands at another fold bifurcation. This causes upstream B-tipping which creates a spike in the downstream response. We observe again in Figure~\ref{fig:2par_bif_diag_loc}(i) that the brown dashed line does not intersect at the tangential meeting of the fold bifurcation curves, but for a slightly larger coupling strength $b$. Increasing the coupling strength moves the two lower fold bifurcations further apart. Between Figures~\ref{fig:2par_bif_diag_loc}(g) and (h), the fold bifurcation to the right moves beyond $\lambda=4$; this corresponds to the end of an overshoot of the downstream threshold. The fold bifurcation retreating to the left has reached a sufficiently low value of $\lambda$ in Figure~\ref{fig:2par_bif_diag_loc}(h), such that downstream B-tipping occurs before upstream B-tipping. After downstream B-tipping, the upstream systems reaches an intermediate steady state that is not accessible for the weaker coupling strengths considered. The subsequent upstream B-tipping is subsequently initiated by crossing another fold bifurcation as the external forcing continues to increase.   

\subsection{Tipping regimes for localised coupling}

\begin{figure}[!ht]
    \centering
    \includegraphics[width=\textwidth]{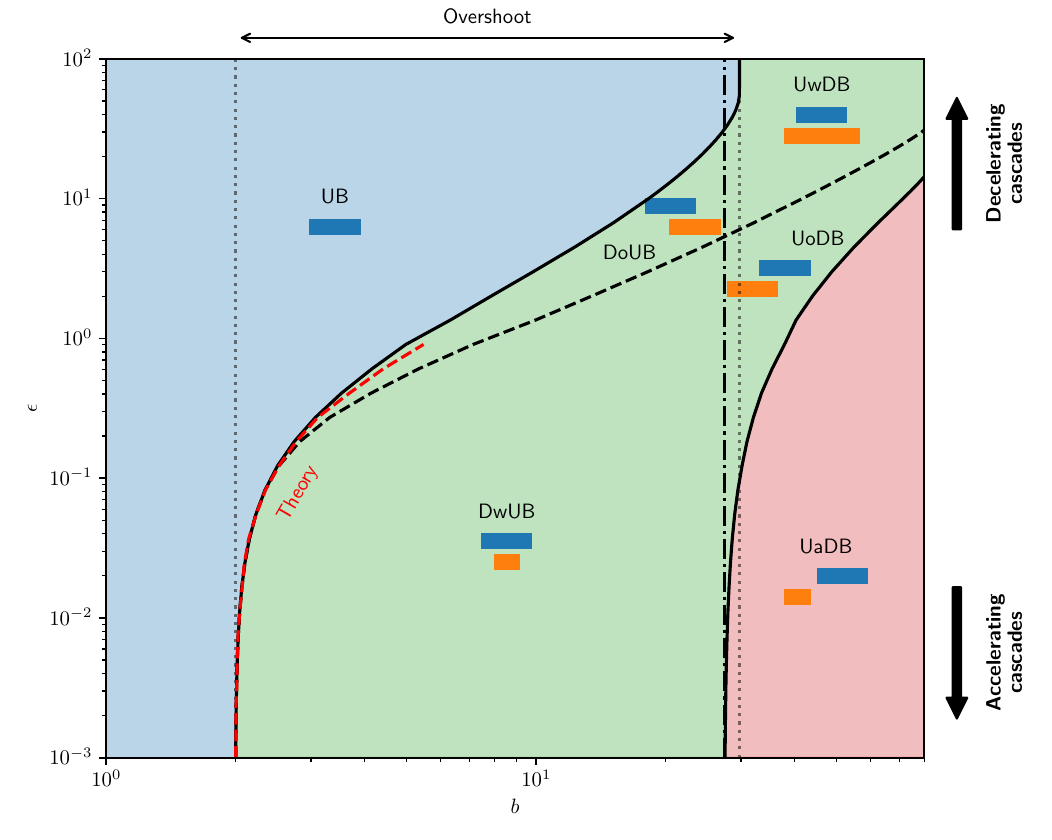}
    
    \caption{\textbf{Tipping regimes for localised coupling} Tipping regimes for system~\eqref{eq:model} with localised coupling $M_2(x)$ dependent on the coupling strength, $b$, and the timescale separation, $\epsilon$. A temporary overshoot of the downstream threshold occurs between the two vertical dotted lines. For increasing $\epsilon$ this overshoot can be progressively further, with asymptotic behaviour shown by the red line marked ``Theory'' derived in \cite{ritchie2019inverse}. 
    The labelling is otherwise as in Figure~\ref{fig:tipping_regimes_lin}.		
    %Blue region (UB): upstream B-tipping only (downstream tracking). Green region represents no intermediate state reach and is partitioned by the relative timings of the tipping onset (vertical dashed line) and offset (curved dashed line) for the upstream and downstream systems. This gives the regions (DwUB): downstream within upstream B-tipping; (DoUB): downstream overlap upstream B-tipping; (UoDB): upstream overlap downstream B-tipping; (UwDB): upstream within downstream B-tipping. Red region (UaDB): upstream after downstream B-tipping that reaches an intermediate state. The boundaries are computed as outlined in Section~\ref{sec:depend}. %Coloured blocks provide schematic illustrations for the tipping duration of the upstream (blue) and downstream (orange) systems. 
    Parameter values as given in Table~\ref{tab:cases} except for the tracking boundary which is computed using $w=1$.
    }
    \label{fig:tipping_regimes_loc}
\end{figure}

The different tipping regimes, as in Figure~\ref{fig:tipping_seq}, for localised coupling are shown in Figure~\ref{fig:tipping_regimes_loc} depending on the timescale separation $\epsilon$ and the coupling strength $b$. These are numerically computed as before, i.e. as outlined in Appendix~\ref{app:boundaries}. If the coupling strength is sufficiently weak ($b<2$), the downstream threshold is not crossed, and UB-tipping (labelled UB; blue region), corresponding to upstream B-tipping with downstream tracking, is the only possible regime. 

For stronger coupling strengths (between the two vertical dotted lines), there is a temporary overshoot of the downstream threshold. For an accelerating cascade, DwUB-tipping predominantly occurs, with downstream B-tipping taking place within upstream B-tipping. In contrast, for decelerating cascades, downstream overshoot with tracking becomes more plausible, corresponding to UB-tipping (the blue region extends to the right for large $\epsilon$). 
Indeed, for small $r$ and $\epsilon$ we can confirm that the hypotheses of Theorem~\ref{thm:criteria} case 2 hold, implying that DwUB-tipping occurs.
Note that the blue region in Figure~\ref{fig:tipping_regimes_loc} can only extend to the right while an overshoot is possible (up to the right vertical dotted line). The boundary that separates the blue and green regions is well approximated by the inverse square law for overshoots, which specifies how long and how far a threshold can be exceeded before tipping occurs \cite{ritchie2019inverse}.

For localized coupling functions such that the coupling is active only during upstream B-tipping, the threshold for the downstream system can never be crossed after the upstream system has completed tipping. This is the case for the coupling $M_2(x)$ and parameters considered here. Therefore, DaUB-tipping, present for linear coupling, is absent for localised coupling. However, similar to the linear coupling case, the green region is still divided into four by the boundaries for the respective timings of the onset (vertical black dash-dot line) and offset (curved black dashed line) for the upstream and downstream systems. The region corresponding to UaDB-tipping is similar to that for linear coupling because there is no coupling from the downstream system to the upstream system.

\section{Discussion}
\label{sec:discuss}

Previous studies have considered the possibility of one element tipping causing another element to tip (and potentially more), in a process commonly known as a tipping cascade \cite{wunderling2024climate,Klose2021Cascade,eker2024cross}. Using a conceptual model of two coupled hysteresis elements, this paper highlights the importance of timescales and coupling for classifying scenarios for a tipping cascade. 
In this study, we considered a single external forcing profile with a fixed timescale separation between the forcing and the upstream system. Despite a much slower forcing timescale than the upstream dynamics, the upstream B-tipping is not instantaneous. Therefore, tipping of one system (if it occurs) may come before or after the other system has tipped and may overlap (or even be contained within) the tipping of the other system. On the other hand, for a decelerating cascade of tipping elements, the upstream B-tipping may be completely contained within the downstream B-tipping. 

If the coupling between the elements is sufficiently weak, then downstream B-tipping will be avoided. However, for very strong coupling and/or if the downstream system is already close to a threshold, downstream B-tipping can occur before upstream B-tipping. 
%Due to our simplifying assumptions of forcing only entering the upstream system and no back-coupling from the downstream to the upstream system, this means that causality approaches \cite{runge2018causal} will not be applicable to distinguish between an upstream system apparently influencing a downstream system (DaUB) but a downstream system not causing B-tipping of an upstream system (UaDB). 

We expect that including a small amount of back coupling and/or time-dependent forcing of the downstream system will not change the results or classification qualitatively, except that degenerate bifurcations can be unfolded. For larger amounts of back coupling, much richer behaviour will be possible, including oscillatory behaviour in the coupled system.

The coupling between two elements can take many forms. Here, we considered both a linear coupling and a localised coupling, which in the latter case can lead to an overshoot of the downstream threshold. If the timescale of the downstream system is slow compared to the upstream system (i.e. $\epsilon \gg 1$ -- a decelerating cascade of tipping elements), then downstream B-tipping can more easily be avoided than for an accelerating cascade of tipping elements ($\epsilon \ll 1$). Furthermore, a localised coupling leads to extrapolation problems for using early warning signals to predict downstream B-tipping, especially in the case of DwUB-tipping for an accelerating cascade of tipping elements \cite{ashwin2025early}. Note that DwUB-tipping is also relevant for the simpler case of linear coupling.

In this study, we considered the relatively simple setup of a single element coupled unidirectionally to another, which still reveals multiple subtleties in the different tipping scenarios that can arise. Though more complex tipping classifications could be unveiled by considering bidirectional coupling and/or coupling additional elements. Alternatively, considering the scenario of a strongly forced upstream system with overshoot in an accelerating cascade of tipping elements may reveal further novel tipping mechanisms. Allowing elements to tip via different tipping mechanisms, such as rate-induced tipping, would likely add a further layer of complexity. 

The roles of timescales and coupling type are important in understanding the propagation of tipping in applications. For instance, the Atlantic Meridional Overturning Circulation (AMOC) is often viewed as a mediator for many different elements due to its central role in a network of tipping elements \cite{wunderling2024climate}. One link in this network couples the Greenland ice sheet (GIS) to AMOC, mediated by freshwater flux into the North Atlantic from the melting GIS. The flux can be thought of as the time derivative of the ice sheet mass. In such a setting, linear coupling is not appropriate. The flux coupling is rather a time-localised coupling, corresponding to the example summarised in Figure~\ref{fig:tipping_regimes_loc}.
Indeed, a temporary overshoot of an AMOC threshold without tipping is discussed in various contexts; see, for example, \cite{ritchie2024esd}.
The AMOC operates on a faster timescale than the Greenland ice sheet, constituting an accelerating cascade of tipping elements. Therefore, the most likely scenarios, from weak to strong coupling, are UB-, DwUB-, and UaDB-tipping. The precise scenario will depend not only on the coupling strength (which may be difficult to determine quantitatively in relation to Figure~\ref{fig:tipping_regimes_loc}) but also on how far the downstream system is from its bifurcation point. 
Another study of GIS to AMOC coupling \cite{sinet2023amoc} (which additionally involves coupling to a tipping West Antarctic ice sheet (WAIS)) demonstrates that timescales and coupling types play a central role in the tipping behaviour of the AMOC. This study examines possible tipping of the AMOC through competition between stabilizing coupling during tipping of the WAIS and destabilizing coupling during tipping of the GIS.

Tipping cascades may be desirable if a chain of tipping transitions is necessary to create a positive change \cite{eker2024cross,nijsse2025leverage}. In such cases, we suggest that timescales and coupling types will be equally important in understanding when a cascade of systems needs to tip for a positive change to be realised downstream.
In summary, tipping is never instantaneous, and systems in cascades can potentially interact in many ways. Hence, a discussion of timescales and coupling types is central to making informed decisions for (and understanding the risks inherent in) coupled systems with the potential to tip.

~

\subsection*{Acknowledgements}

We thank the reviewers for their detailed and helpful comments and suggestions. This research was supported by the European Union's Horizon Europe research and innovation programme under grant agreements No. 101137601 (ClimTip) and and 101184070 (P2F).
We acknowledge support by the UK Advanced Research and Invention Agency (ARIA) via the project ``AdvanTip". A.S.vdH acknowledges funding by the Dutch Research Council (NWO) through the NWO-Vici project ``Interacting climate tipping elements: When does tipping cause tipping?'' (project VI.C.202.081). Views and opinions expressed are however those of the author(s) only and do not necessarily reflect those of the European Union or the European Climate, Infrastructure and Environment Executive Agency (CINEA) or other funders. Neither the European Union nor the granting authority can be held responsible for them. 
For the purpose of open access, the authors have applied a Creative Commons Attribution (CC BY) licence to any Author Accepted Manuscript version arising from this submission.

\subsection*{Data accessibility statement}
The code to generate the data for the figures presented is available on
\url{https://github.com/PaulRitchie89/tippingcascades_timescales}. 

\bibliographystyle{plain}

\bibliography{cascadebib}

\appendix

\section{Numerical classification of boundaries between tipping scenarios}
\label{app:boundaries}

If both upstream and downstream systems undergo tipping, we can numerically approximate boundaries between the scenarios of Figure~\ref{fig:tipping_seq} for the system (\ref{eq:model}). We do this by finding where various combinations of onset and offset for the upstream and downstream are simultaneous. For example, if there is a $t_1$ such that $\Lambda(rt_1)=M(x(t_1))=\lambda_u$, this corresponds to the boundary of simultaneous onset of tipping, while if $t_1$ is such that $\Lambda(rt_1)=\lambda_u$ and $y(t_1) = w$, this corresponds to the boundary where upstream tipping starts simultaneously with downstream tipping ending. The other boundaries are found in a similar manner. The upstream dynamics are the same throughout, and so with the upstream timings identified, a bisection method is used to align the downstream timing. These boundary locations clearly depend on the choice of threshold. The code is available (see data accessibility statement).

\begin{figure}[!ht]
    \centering
    \includegraphics[width=0.8\textwidth]{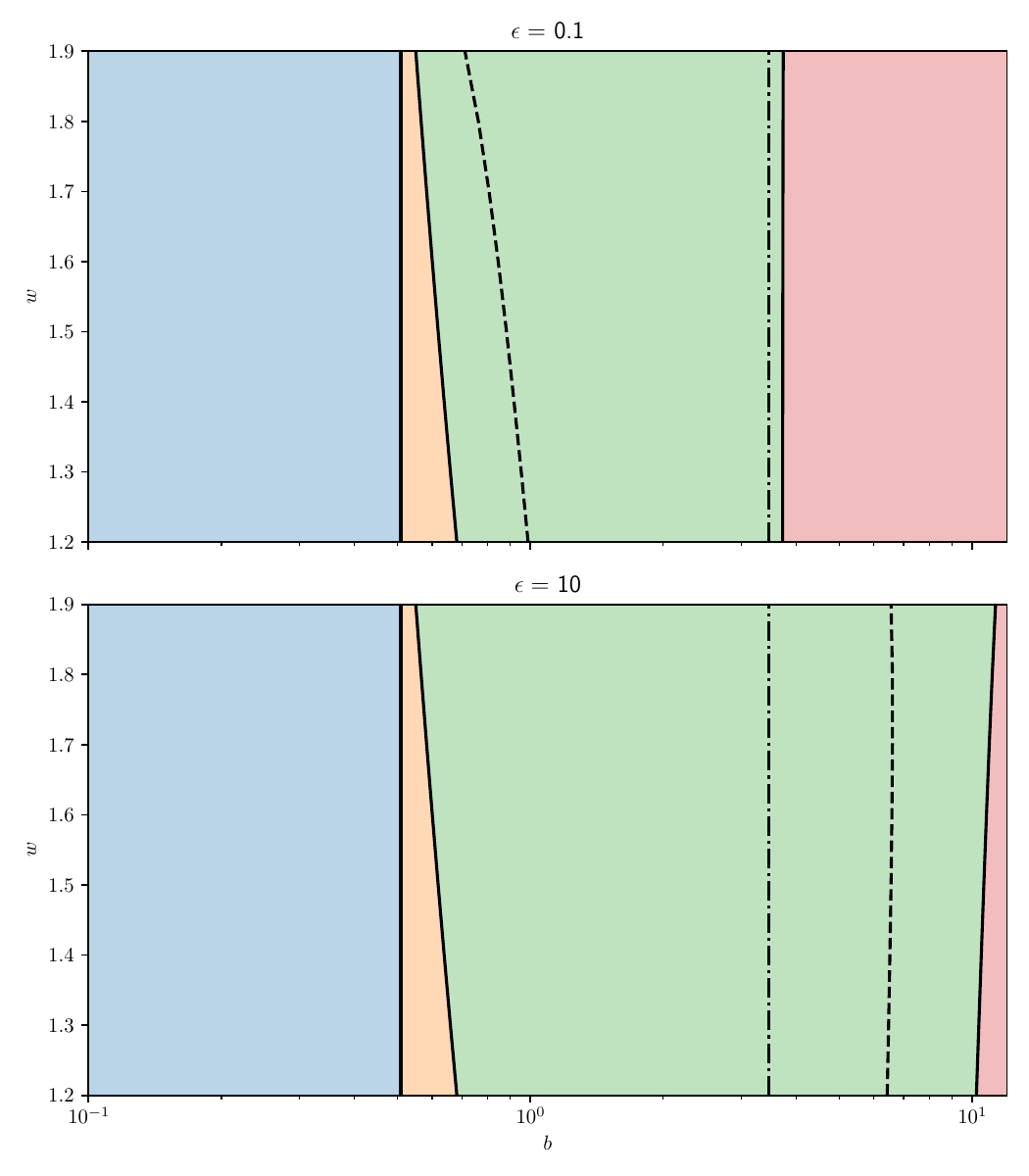}
    
    \caption{\textbf{Sensitivity analysis of tipping boundaries based on the offset threshold for linear coupling.} Tipping boundaries for system~\eqref{eq:model} with linear coupling $M_1(x)$ dependent on the coupling strength $b$, and the offset threshold $w$ for an accelerating cascade (top, $\epsilon=0.1$) and a decelerating cascade (bottom, $\epsilon=10$) of tipping elements. Blue region: upstream B-tipping only (downstream tracking). Orange region: downstream B-tipping after upstream B-tipping. The green region represents when there is at least a brief period of both upstream and downstream systems tipping. Simultaneous tipping onset (resp. offset) for the upstream and downstream systems are given by dashed-dot (resp. dashed) line. Red region: upstream after downstream B-tipping that reaches an intermediate state. The boundaries are computed as outlined in Section~\ref{sec:depend}. Parameter values as given in Table~\ref{tab:cases}.
    }
    \label{fig:tipping_regimes_lin_wsens}
\end{figure}

Note that the classification will become independent of choice of threshold only for limiting timescale ratios; for the finite timescale separations considered here there is nonetheless only a weak dependence on choice of threshold: see Figure~\ref{fig:tipping_regimes_lin_wsens} which illustrates the dependence of the boundaries in Figure~\ref{fig:tipping_regimes_lin} on the threshold $w$ for cases $\epsilon=0.1$ and $\epsilon=10$.

\end{document}